# Surgery and duality

By Matthias Kreck

## 1. Introduction

Surgery, as developed by Browder, Kervaire, Milnor, Novikov, Sullivan, Wall and others is a method for comparing homotopy types of topological spaces with diffeomorphism or homeomorphism types of manifolds of dimension $\geq 5$. In this paper, a modification of this theory is presented, where instead of fixing a homotopy type one considers a weaker information. Roughly speaking, one compares n-dimensional compact manifolds with topological spaces whose $k$-skeletons are fixed, where $k$ is at least $[n/2]$. A particularly attractive example which illustrates the concept is given by complete intersections. By the Lefschetz hyperplane theorem, a complete intersection of complex dimension $n$ has the same $n$-skeleton as $\mathbb{CP}^n$ and one can use the modified theory to obtain information about their diffeomorphism type although the homotopy classification is not known. The theory reduces this classification result to the determination of complete intersections in a certain bordism group. This was under certain restrictions carried out in [Tr]. The restrictions are: If $d = d_1 \cdot \ldots \cdot d_r$ is the total degree of a complete intersection $X_{d_1,\ldots,d_r}^n$ of complex dimension $n$, then the assumption is, that for all primes $p$ with $p(p-1) \leq n+1$, the total degree $d$ is divisible by $p^{[(2n+1)/(2p-1)]+1}$.

THEOREM A.   *Two complete intersections $X_{d_1,\ldots,d_r}^n$ and $X_{d_1',\ldots,d_s'}^n$ of complex dimension $n > 2$ fulfilling the assumption above for the total degree are diffeomorphic if and only if the total degrees, the Pontrjagin classes and the Euler characteristics agree.*

Note that the $k$th Pontrjagin class is a multiple of $x^{2k}$, where $x$ generates the second cohomology of the complete intersection. Thus we can compare this invariant for different complete intersections. There are explicit formulas for all these invariants. It is open whether this theorem holds for arbitrary complete intersections of complex dimension $> 2$.

The $k$-skeleton is not an invariant of a topological space and thus we pass to the closely related language of Postnikov towers. The *normal $k$-type* of a manifold is the fibre homotopy type of a fibration $B \to BO$ such that $\pi_i(B \to BO) = 0$ for $i \geq k + 2$, admitting a lift of the normal Gauss map $\nu : M \to BO$ to a map $\bar\nu : M \to B$ such that $\pi_j(\bar\nu : M \to B) = 0$ for $j \leq k+1$. We call such a lift a *normal k-smoothing*. A normal $k$-smoothing determines an element in an obvious bordism group given by the normal $k$-type. The main result of this paper concerning the classification of manifolds is the following.

THEOREM B. *Let $k \geq [n/2] - 1$. A normal bordism $W$ of dimension $n + 1 > 4$ between two normal k-smoothings on manifolds $M_0$ and $M_1$ with the same Euler characteristics is bordant to an s-cobordism if and only if an algebraic obstruction $\theta(W)$ is elementary. Thus $M_0$ and $M_1$ are diffeomorphic, if $n > 4$, and homeomorphic, if $n = 4$ and $\pi_1$ is good in the sense of* [Fr2].

In the most general case the obstruction $\theta(W)$ lies in a monoid depending only on the fundamental group and the orientation character given by the first Stiefel-Whitney class. If $k \geq [n/2]$ the obstruction is contained in a subgroup of the monoid and one obtains as special cases the Wall-obstructions and classification results. For a detailed formulation of Theorem B we refer to Theorem 3 (Section 5) and Theorem 4 (Section 6). For simply connected manifolds a similar approach to the classification problem was carried out by M. Freedman [Fr1].

The obstructions are particularly complicated in the extreme case $k = [n/2] - 1$, even for simply connected manifolds. It is surprising that they can be omitted if the manifolds are of dimension $2q$ and one allows stabilization by $S^q \times S^q$. Then the result is the following which generalizes a result by Freedman for specific 1-connected manifolds [Fr1, Th. 3].

THEOREM C. *Two closed 2q-dimensional manifolds with the same Euler characteristic and the same normal $(q - 1)$-type, admitting bordant normal $(q - 1)$-smoothings, are diffeomorphic after connected sum with $r$ copies of $S^q \times S^q$ for some $r$.*

If the fundamental group is finite, one has cancellation results. In joint work with Ian Hambleton [H-K1, Th. 1.3] we showed that up to homeomorphism one can take $r = 2$ for $q = 2$ and a similar argument holds up to diffeomorphism for $q > 2$. If $q = 2$, the main theorem of [H-K3] gives cancellation up to homeomorphism down to $r = 1$.

In Section 7, Corollary 4, we will give a short proof of the cancellation result for $q > 2$ by using the unitary stability techniques of Bass [Ba2] (as in [H-K1]) to analyse the obstruction $\theta(W)$ directly. For 1-connected manifolds, elementary arguments give the following result which is best possible, and which for $q$ odd was proved by Freedman [Fr1, Th. 1]:

THEOREM D. *For $q > 2$, two closed simply connected $2q$-dimensional manifolds $M_0$ and $M_1$ with the same Euler characteristic and the same normal $(q-1)$-type admitting bordant normal $(q-1)$-smoothings are diffeomorphic if either $q$ is odd or $q$ is even and $M_0 = M_0'\sharp S^q \times S^q$.*

This is the background for Theorem A. To mention an application of the cancellation results to nonsimply connected manifolds we combine them with the exact surgery sequence [W1] to compute under certain assumptions the group of connected components of (local) orientation-preserving simple homotopy self equivalences $\pi_0(\mathrm{Aut}(M))$ modulo the group $\pi_0(\mathrm{Diff}(M))$ of isotopy classes of (local) orientation-preserving diffeomorphisms in terms of an exact sequence.

THEOREM E. *Let $M^{2q}$ either be 1-connected with $q$ odd or 1-connected with $q$ even and $M = M'\sharp S^q \times S^q$, or $\pi_1(M)$ finite and $M = M'\sharp 2(S^q \times S^q)$. If $q > 2$ there is an exact sequence*

$$[\Sigma(M), G/O] \to L_{2q+1}^s(\pi_1(M), w_1(M))$$
$$\to \pi_0(\mathrm{Aut}(M))/\pi_0(\mathrm{Diff}(M)) \to [M, G/O] .$$

To give an application of Theorem B to manifolds with infinite fundamental groups we present a very quick proof of the following result, which was independently proven by Freedman and Quinn [F-Q, Th. 10.7.A].

THEOREM F. *Two closed topological 4-dimensional spin 4-manifolds with infinite cyclic fundamental group are homeomorphic if and only if they have isometric intersection forms on $\pi_2$.*

The intersection form is a quadratic form with values in $\Lambda = \mathbb{Z}[\pi_1]$, which is described in Section 5.

For odd-dimensional manifolds I do not know a result like Theorem C and so it is necessary to analyse the obstructions $\theta(W)$. We will carry this out in a special case which has applications to the classification of 1-connected 7-dimensional homogeneous spaces. These homogeneous spaces have torsion-free second homology group and isomorphic finite fourth cohomology group generated by the first Pontrjagin class and decomposable classes. The normal 2-type is then determined by the second Betti number and the second Stiefel-Whitney class. An analysis of $\theta(W)$ for a bordism between two such homogeneous spaces leads to:

THEOREM G. *Two 1-connected 7-dimensional homogeneous spaces with the cohomolgical properties above are diffeomorphic if and only if they have equal second Betti number and Stiefel-Whitney class and if there is a bordism $W$ between normal 2-smoothings such that $\text{sign}(W) = 0$ and the characteristic numbers $p_1^2(W)$, $z_i^2 p_1(W)$ and $z_i^2 z_j^2$ vanish, where $z_i$ are classes in $H^2(W; \mathbb{Z})$ restricting to a basis of each of the two boundary components.*

In joint work with Stephan Stolz we analysed this situation further and showed that besides the second Betti number and Stiefel-Whitney class certain spectral invariants determine the diffeomorphism type. An explicit calculation of the spectral invariants gave the first examples of homeomorphic but not diffeomorphic homogeneous spaces [Kr-St1] and of manifolds where the moduli space of metrics with positive sectional curvature is not connected [Kr-St2].

Besides the aim of obtaining explicit classification results the theory sheds some light on the role of duality for manifolds. Poincaré duality reflects some symmetry between the $k$ and $n-k$ handles of a compact manifold. Prescribing the $[n/2]$-skeleton and classifying the corresponding manifolds shows how far manifolds are determined by their handles up to half the dimension. We will mention a result which demonstrates that in a particular situation even the cohomology ring up to the middle dimension plus the Pontrjagin classes and a certain homology class determine the manifolds up to finite ambiguity. Sullivan [Su] introduced minimal models and the notion of a formal space, which means that the minimal model is determined by its rational cohomology ring. We abbreviate for an $n$-dimensional manifold $M$ the truncated cohomology ring $\sum_{i \leq [n/2]+1} H^i(M; \mathbb{Z})$ by $H_{\leq [n/2]+1}(M)$ and the subalgebra of the real minimal model of $H_{\leq [n/2]+1}(M)$ generated by elements of degree $\leq [n/2]$ by $\mathcal{M}_{[n/2]}(M)$. The fundamental class of $M$ determines a class $\alpha(M) \in H_n(\mathcal{M}_{[n/2]}(M))$. The result which was proved in [Kr-Tr] using the modified surgery theory is the following. Let $n \geq 5$. The diffeomorphism type of a 1-connected closed smooth $n$-manifold with formal $([n/2]+1)$-skeleton is determined up to finite ambiguity by the truncated cohomology ring $H_{\leq [n/2]+1}(M)$, the real Pontrjagin classes and the class $\alpha(M) \in H_n(\mathcal{M}_{[n/2]}(M))$.

Most of the results of this paper were obtained in the early eighties and were circulated as [Kr3]. A plan to write a monograph based on this preprint could not yet be realized. Since the theory was meanwhile used in several papers ([Be], [Da], [F-K], [F-K-V], [H-K1], [H-K2], [H-K3], [H-K4], [H-K-T], [K-L-T], [Kr4], [Kr-St1], [Kr-St2], [Kr-St3], [Kr-Tr], [Sto2], [Te], [Tr], [Wa]), I decided to publish the most important results in the present form.

I would like to thank Stephan Stolz and Peter Teichner for many helpful discussions about the theory, and the referee for detailed suggestions improving the presentation.

## 2. Normal $k$-smoothings

We will formulate and prove our general results in the smooth category. Most results can with appropriate modifications be proved in the piecewise linear or topological category. (Replace the differential normal bundle by the corresponding PL - or TOP bundle.) This follows from the basic results of [K-S].

We use the language of manifolds with $B$-structures. Here $B$ is a fibration over $BO$ and a *normal $B$-structure* on an $n$-dimensional manifold $M$ in $B$ is a lift $\bar{\nu}$ of the stable normal Gauss map $\nu : M \to BO$ to $B$. Since the normal Gauss map depends on an embedding of $M$ into $\mathbb{R}^{n+r}$ for $r$ large, one has to interpret this with care and we refer to [St, p. 14 ff] for details. Since we will frequently use homotopy groups we equip all spaces, without special mentioning, with base points and assume that maps preserve the base points. In particular if we orient the classifying bundle over $BO$ at the base point the normal Gauss map induces a local orientation at the base point and so all orientable connected manifolds come with a given orientation.

*Definition.* Let $B$ be a fibration over $BO$.
i) A normal $B$-structure $\bar{\nu} : M \to B$ of a manifold $M$ in $B$ is a *normal $k$-smoothing*, if it is a $(k+1)$-equivalence.
ii) We say that $B$ is $k$-*universal* if the fibre of the map $B \to BO$ is connected and its homotopy groups vanish in dimension $\geq k+1$.

Obstruction theory implies that if $B$ and $B'$ are both $k$-universal and admit a normal $k$-smoothing of the same manifold $M$, then the two fibrations are fibre homotopy equivalent. Furthermore, the theory of Moore-Postnikov decompositions implies that for each manifold $M$ there is a $k$-universal fibration $B^k$ over $BO$ admitting a normal $k$-smoothing of $M$. For background on these basic homotopy theoretic facts we refer to [Ste] or more generally to [Bau]. Thus the fibre homotopy type of the fibration $B^k$ over $BO$ is an invariant of the manifold $M$ and we call it the *normal $k$-type* of $M$ denoted $B^k(M)$. We note that if two manifolds have homotopy equivalent $(k+1)$-skeletons and isomorphic normal bundles over them, then they have the same normal $k$-type. By obstruction theory one obtains a classification of all normal $k$-smoothings of $M$ in $B^k(M)$. The group of fibre homotopy classes of fibre homotopy self-equivalences $\mathrm{Aut}(B^k(M))$ acts effectively and transitively on the set of normal $k$-smoothings of $M$.

There is an obvious bordism relation on closed $n$-dimensional manifolds with normal $B$ structures and the corresponding bordism group is denoted $\Omega_n(B)$ [St]. Normal $k$-smoothings give special elements in $\Omega_n(B)$ and these are independent of the choice of the normal $k$-smoothing in $\Omega_n(B)/\mathrm{Aut}(B)$.

*Remark.* If $k$ is larger than $n$, the dimension of $M$, then $B^k(M)$ is equivalent to the normal homotopy type of $M$: Two manifolds have the same (= fibre homotopy equivalent) normal $k$-type if and only if there is a homotopy equivalence preserving the normal bundle. Thus the starting point of the original surgery theory, the normal homotopy type, is a special case of our setting.

We will demonstrate now using some examples that it is often much easier to determine the normal $[n/2]-1$-type of a manifold than its normal homotopy type.

Consider an $n$-dimensional homotopy sphere $\Sigma$. To describe the normal $k$-type of $\Sigma$ we need the following notion. Let $X$ be a connected CW-complex. The *$k$-connected cover $X\langle k \rangle$* is a CW-complex which up to homotopy equivalence is characterized by the property that $X\langle k \rangle$ is $k$-connected and there is a fibration $p : X\langle k \rangle \longrightarrow X$ inducing isomorphisms on $\pi_i$ for $i > k$.

PROPOSITION 1. *Let $\Sigma^n$ be an $n$-dimensional homotopy sphere and $k < n - 1$. Then the normal $k$-type of $\Sigma^n$ is the fibration $p : BO\langle k+1 \rangle \longrightarrow BO$.*

*Proof.* Since the fibration $p : BO\langle k+1 \rangle \longrightarrow BO$ induces an isomorphism on $\pi_n$ for $n > k$, the normal Gauss map lifts and the lift is automatically a $(k+1)$-equivalence. □

*Remark.* For $k \geq n$ the normal $k$-type is equivalent to the normal homotopy type of a homotopy sphere. The determination of this is an important step in the analysis of homotopy spheres by ordinary surgery theory as was done by Kervaire and Milnor [K-M]. The additional information needed for this is that the stable normal bundle of a homotopy sphere is trivial [K-M, Th. 3.1]. It should be noted that the proof of this fact is not elementary (it uses the Hirzebruch signature theorem as well as Adams's result about the injectivity of the $J$-homomorphism and of course Bott periodicity). In contrast, the proof of the proposition for $k < n - 1$ is completely elementary. One can, based on this completely elementary proposition, see that one gets the same information about the diffeomorphism classification of homotopy spheres as Kervaire and Milnor.

Next, we determine the normal 1-type of a compact manifold. This is relevant for determining the homeomorphism type of compact 4-manifolds and for applications to manifolds of dimension $> 4$ with metric of positive scalar curvature. Consider triples $(\pi, w_1, w_2)$ where $\pi$ is a finitely presentable group and $w_i \in H^i(K(\pi, 1); \mathbb{Z}/2)$ are cohomology classes. Two such triples are called *isomorphic* if there is an isomorphism $f : \pi \longrightarrow \pi'$ such that $f^* w_i' = w_i$. We denote the isomorphism class by $[\pi, w_1, w_2]$. Similarly we introduce isomorphism classes of pairs $[\pi, w_1]$, where $w_1$ is an element of $H^1(K(\pi, 1); \mathbb{Z}/2)$.

Given $(\pi, w_1)$ we consider the real line bundle $E \to K(\pi, 1)$ with $w_1(E) = w_1$. Consider the composition

$$K(\pi, 1) \times BSO \xrightarrow{E \times p} BO \times BO \xrightarrow{\oplus} BO,$$

where $E: K(\pi, 1) \to BO$ is the classifying map of the stable bundle given by $E$ and $\oplus$ is the $H$-space structure on $BO$ given by the Whitney sum. We denote the corresponding fibration by $B[\pi, w_1]$. The normal Gauss map $\nu: M \to BO$ together with $u: M \to K(\pi, 1)$ determines a lift $\bar{\nu}: M \to B(\pi, w_1)$ of $\nu$ and it is easy to check that $\bar{\nu}$ is a 2-equivalence.

Given $(\pi, w_1, w_2)$, we consider the following pullback square

$$B(\pi, w_1, w_2)@>>> K(\pi, 1)@.$$

$$@VpVV @VVw_1 \times w_2 V@.$$

$$BO@>>w_1(EO) \times w_2(EO)> K(\mathbb{Z}/2, 1) \times K(\mathbb{Z}/2, 2),$$

where $w_i(EO)$ are the Stiefel-Whitney classes of the universal bundle. The fibre homotopy type of $p: B(\pi, w_1, w_2) \longrightarrow BO$ is determined by the isomorphism class of $(\pi, w_1, w_2)$ and is denoted by $B[\pi, w_1, w_2]$.

If $M$ is a compact manifold (implying $\pi_1(M)$ is finitely presentable) and $u: \pi_1(M) \to \pi$ is an isomorphism we denote the corresponding map $M \to K(\pi, 1)$ again by $u$ ($u$ is unique up to homotopy and a classifying map of the universal covering). If $w_2(\widetilde{M}) = w_2(\nu(\widetilde{M})) = 0$ there are unique classes $w_i \in H^i(K(\pi, 1); \mathbb{Z}/2)$ with $u^*w_i = w_i(\nu(M))$ for $i = 1, 2$. This is clear for $i = 1$ and for $i = 2$ one uses the short exact sequence [Bro]:

$$0 \to H^2(K(\pi, 1); \mathbb{Z}/2)@>u^*>> H^2(M; \mathbb{Z}/2)@>p^*>> H^2(\widetilde{M}; \mathbb{Z}/2).$$

Obviously $[\pi, w_1, w_2]$ is an invariant of $M$.

The normal Gauss map $\nu: M \to \mathcal{B}$ together with $u: M \to K(\pi, 1)$ determines a lift $\bar{\nu}: M \to B(\pi, w_1, w_2)$ of $\nu$ and it is easy to check that $\bar{\nu}$ is a 2-equivalence. We summarize these considerations as:

PROPOSITION 2. *If* $w_2(\widetilde{M}) \neq 0$ *then the normal 1-type of a compact manifold* $M$ *is* $B[\pi, w_1]$, *and if* $w_2(\widetilde{M}) = 0$ *then it equals* $B[\pi, w_1, w_2]$.

Finally we determine the normal $(n - 1)$-type of a complete intersection. Let $f_1, \ldots, f_r$ be homogeneous polynomials on $\mathbb{C}P^{n+r}$ of degree $d_1, \ldots, d_r$. If the gradients of these polynomials are linearly independent, the set of common zeros is a smooth complex manifold of complex dimension $n$, a nonsingular complete intersection. As was noted by Thom, the diffeomorphism type of nonsingular complete intersections depends only on the unordered tuple $(d_1, \ldots, d_r)$ called the multi-degree. We denote this diffeomorphism type by $X^n_{d_1, \ldots, d_r}$. It is natural to ask for a diffeomorphism classification of this very interesting class of algebraic manifolds.

Except under some restrictive assumptions [L-W1], [L-W2], even the homotopy classification of the $X_{d_1,\ldots,d_r}^n$'s is unknown, which is the first step in the ordinary surgery theory. On the other hand the topology of $X_n(d)$ up to half the dimension is known. According to Lefschetz the inclusion

$$i : X_{d_1,\ldots,d_r}^n \longrightarrow \mathbb{C}P^\infty$$

is an $n$-equivalence.

Moreover, it is easy to see that the normal bundle of $X_{d_1,\ldots,d_r}^n$ is isomorphic to

$$\nu(X_{d_1,\ldots,d_r}^n) \cong i^*(\nu(\mathbb{C}P^{n+r}) \oplus H^{d_1} \oplus \cdots \oplus H^{d_r})$$
$$\cong i^*(-(n+r+1) \cdot H \oplus H^{d_1} \oplus \cdots \oplus H^{d_r})$$

where $H$ is the Hopf bundle and $H^{d_i}$ means the $d_i$-fold tensor product. We abbreviate $\delta = (d_1, \ldots d_r)$. Denote the classifying map of

$$-(n+r+1) \cdot H \oplus H^{d_1} \oplus \cdots \oplus H^{d_r}$$

by $\xi(n,\delta) : \mathbb{C}\mathbb{P}^\infty \to BO$. We transform the composition of $\xi(n,\delta) \times p$ $: \mathbb{C}P^\infty \times BO\langle n+1 \rangle \to BO \times BO$ and the Whitney sum $\oplus : BO \times BO \to BO$ into a fibration and denote the projection map of this fibration by $\xi(n,\delta) \oplus p$ $: \mathbb{C}P^\infty \times BO\langle n+1 \rangle \to BO$. Then by construction the normal Gauss map of $X_{d_1,\ldots,d_r}^n$ admits a lift over this fibration by a $n$-equivalence. Then:

PROPOSITION 3. *The normal $(n-1)$-type of a complete intersection $X_n(\delta)$ is*

$$\mathbb{C}P^\infty \times BO\langle n+1 \rangle @>\xi(n,\delta) \oplus p>> BO.$$

## 3. Surgery below the middle dimension and first applications

In homotopy theory one can, for a topological space $X$ and $r \geq 1$, eliminate arbitrary elements $[f] \in \pi_r(X)$ by attaching an $(r+1)$-cell via $f$. More precisely consider $Y = D^{r+1} \cup_f X$. Then the inclusion $i : X \to Y$ is an $r$-equivalence and $[f]$ with all its translates under the action of $\pi_1(X)$ generates the kernel of $i_* : \pi_r(X) \to \pi_r(Y)$ (for $r > 1$ see [Wh], for $r = 1$ this follows from van Kampen's theorem).

Surgery is an attempt to do constructions which have the same effect on homotopy groups within the category of manifolds [Br], [W1]. To stay within the category of manifolds, we start with an embedding

$$f : S^r \times D^{m-r} \hookrightarrow \mathring{M}$$

where $M$ is a $m$-dimensional manifold. Then we define

$$W := D^{r+1} \times D^{m-r} \underset{f}{\cup} M \times I$$

where we consider $f$ as a map to $M \times \{1\}$.

$W$ is a manifold with corners but we will always straighten the angles occurring at $f(S^r \times S^{m-r-1})$ [C-F]. This construction is called *attaching an* $(r+1)$-*handle* and $W$ the *trace of a surgery via* $f$.

The boundary of $W$ is $M \cup (\partial M \times I) \cup M'$ and we call $M'$ the *result of a surgery of index* $r + 1$ *via* $f$. More explicitly,

$$M' = D^{r+1} \times S^{m-r-1} \underset{f}{\cup} (M - f(S^r \times \overset{\circ}{D}{}^{m-r})).$$

Obviously $W$ is homotopy equivalent to $Y = D^{r+1} \underset{f|_{S^r \times \{0\}}}{\cup} M$, the result of attaching a cell via $f|_{S^r \times \{0\}}$. From the construction of $W$ and $M'$ it is not difficult to see that $W$ can also be viewed as the trace of a surgery on $M'$ via the obvious embedding of $D^{r+1} \times S^{m-r-1}$ into $M'$ [Mi1]. In particular, $W$ is homotopy equivalent to $Y' = D^{m-r} \underset{\{0\} \times S^{m-r-1}}{\cup} M'$.

The following lemma demonstrates the analogy of the two constructions "attaching a cell" and "surgery" as far as the effect on homotopy groups is concerned.

LEMMA 1. *Let* $f : S^r \times D^{m-r} \hookrightarrow M^m$ *be an embedding into a connected manifold. Let* $W$ *be the trace of a surgery via* $f$ *and* $M'$ *be the result of a surgery via* $f$.

i) *The inclusion* $i : M \to W$ *is an* $r$-*equivalence and* $[f|S^r \times \{0\}]$ *and its translates under the action of* $\pi_1(M)$ *generate the kernel of* $i_* : \pi_r(M) \to \pi_r(W)$.

ii) *The inclusion* $j : M' \to W$ *is an* $(m - r - 1)$-*equivalence and* $[\{0\} \times S^{m-r-1}] \in \pi_{m-r-1}(M')$ *and its translates under the action of* $\pi_1(M')$ *generate the kernel of* $j : \pi_{m-r-1}(M') \to \pi_{m-r-1}(W)$.

iii) *If* $k < r$ *and* $k < m - r - 1$, *then*

$$\pi_k(M') \cong \pi_k(M) \cong \pi_k(W)$$

*and, if* $2r < m - 1$

$$\pi_r(M') \cong \pi_r(M)/U$$

*where* $U$ *is generated by* $[f|_{S^r \times \{0\}}]$ *and its translates under the action of* $\pi_1(M)$.

*Proof.* The results follow from [Wh, p. 213] and van Kampen's theorem since

i) $W \simeq D^{r+1} \underset{f|_{S^r \times \{0\}}}{\cup} M$ and

ii) $W \simeq D^{m-r} \underset{\{0\} \times S^{m-r-1}}{\cup} M'$.

iii) follows from i) and ii). ▫

To apply the construction of attaching handles to eliminate elements in $\pi_r(M)$, it is necessary to know which elements in $\pi_r(M)$ can be represented by embeddings $f : S^r \times D^{m-r} \hookrightarrow M$. We have some control over this in the situation described in Section 2. Let

$$\xi : B \to BO$$

be a fibration and $\bar{\nu} : M \to B$ a normal $B$-structure. If $r < \frac{m}{2}$, the Whitney embedding theorem [Hi] implies that any map $S^r \to M$ is homotopic to an embedding $f$. If $[f]$ lies in the kernel of $\bar{\nu} : \pi_r(M) \to \pi_r(B)$, the stable normal bundle of this embedding is trivial. Since the dimension of the normal bundle is greater than $r$, it is actually trivial [Ste]. Thus, we have shown the first part of the following lemma.

LEMMA 2.    *Let $\xi : B \to BO$ be a fibration and $(M, \bar{\nu})$ be a normal $B$-structure.*

i) *If $r < \frac{m}{2}$ any element in the kernel of $\bar{\nu}_* : \pi_r(M) \to \pi_r(B)$ can be represented by an embedding*

$$f : S^r \times D^{m-r} \hookrightarrow M.$$

ii) *Let $f : S^r \times D^{m-r} \hookrightarrow M$ be an embedding representing a homotopy class in the kernel of $\bar{\nu}_\star$. For $1 < r < \frac{m}{2}$, $f$ can be modified by a self-diffeomorphism on $S^r \times D^{m-r}$, so that $\bar{\nu} : M \to B$ extends to a normal $B$-structure of $W$, the trace of the surgery via $f$. Denote the restriction of any such extensions to $M'$, the result of the surgery, by $\bar{\nu}' : M' \to B$.*

iii) *For $1 < r = \frac{m}{2}$ and $r \neq 3, 7$, or $r = 3, 7$ and there is $\beta \in \pi_{r+1}(B)$ with $\beta^* \xi^* w_{r+1} \neq 0$, $w_{r+1} \in H^*(BO; \mathbb{Z}/2)$ the Stiefel-Whitney class, the same statement as in ii) holds.*

*Proof.* We only have to show ii) and iii). The embedding $f : S^r \times D^{m-r} \hookrightarrow M$ induces a normal $B$-structure on $S^r \times D^{m-r}$ denoted by $f^* \bar{\nu}$. There is a unique (up to homotopy) $B$-structure on $D^{r+1} \times D^{m-r}$ and we have to show that, after perhaps modifying the embedding $f$, we can achieve that its restriction to $S^r \times D^{m-r}$ is $f^* \bar{\nu}$. Let $F$ be the fibre of $\xi : B \to BO$. The different $B$-structures on $S^r \times D^{m-r}$ are classified by $\pi_r(F)$, as follows from the long exact homotopy sequence. Since $f^* \bar{\nu}|_{S^r \times \{0\}}$ represents 0 in $\pi_r(B)$ by assumption, the $B$-structures are in the image of the boundary operator $d : \pi_{r+1}(BO) \to \pi_r(F)$. For a map $\alpha : S^r \to O(m-r)$ we consider the diffeomorphism $g_\alpha : S^r \times D^{m-r} \to S^r \times D^{m-r}$ mapping $(x, y) \longmapsto (x, \alpha(x) \cdot y)$. Then $f^* \bar{\nu}$ and $(f \cdot g_\alpha)^* \bar{\nu}$ differ by $d(i\alpha) \in \pi_r(F)$, where $i : O(m-r) \to O$ is the inclusion and we consider $i\alpha$ as an element of $\pi_{r+1}(BO) \cong \pi_r(O)$. Since $m - r > r$, $i_* : \pi_r(O(m-r)) \to \pi_r(O)$ is surjective [Ste] which finishes the proof of ii).

For $r = \frac{m}{2}$ the same argument as above works as long as $i_* : \pi_r(O(r)) \to \pi_r(O)$ is surjective. This is the the case for $r \neq 3, 7$ as follows from results in [Ste]. For $r = 3, 7$ the map is not surjective but has a cokernel $\mathbb{Z}/2$. This cokernel is detected by the Stiefel Whitney class $w_{r+1}$ of the bundle over $S^{r+1}$ classified by an element of $\pi_r(O)$. Looking at the homotopy sequence of the fibration $B \to BO$: $\pi_{r+1}(B) \to \pi_{r+1}(BO) \to \pi_r(F) \to \pi_r(B)$ we see that if there is $\beta \in \pi_{r+1}(B)$ with $\beta^\star \xi^* w_{r+1} \neq 0$ there is no obstruction for finding a diffeomorphism $g_\alpha : S^r \times D^{m-r} \to S^r \times D^{m-r}$, so that after changing the embedding with this diffeomorphism $\bar\nu : M \to B$ extends to a normal $B$-structure of $W$, the trace of the surgery via $f$. ☐

We call an embedding $f : S^r \times D^{m-r} \hookrightarrow M$, where $\bar\nu$ extends to a normal $B$-structure of the trace a *compatible* embedding.

Combining the information about the effect of attaching a cell for homotopy groups with Lemma 1 we get the following result. Before we formulate it recall that the integral group ring $\mathbb{Z}[\pi]$ of a group $\pi$ is the ring of all formal linear combinations $\sum n_g g$, where $g$ runs over elements of $\pi$ and all but finitely many $n_g$ are zero. We abbreviate $\mathbb{Z}[\pi_1(B)]$ by $\Lambda$. If $\pi$ is the fundamental group of a space $X$ then it acts on all homology groups of the universal covering and on all homotopy groups of dimension $> 1$, making these groups into $\Lambda$-modules in such a way that the Hurewicz homomorphism is a $\Lambda$-module homomorphism.

PROPOSITION 4. *Let $\xi : B \to BO$ be a fibration and assume that $B$ is connected and has a finite $[m/2]$-skeleton. Let $\bar\nu : M \to B$ be a normal $B$-structure on an $m$-dimensional compact manifold $M$. Then, if $m \geq 4$, by a finite sequence of surgeries $(M, \bar\nu)$ can be replaced by $(M', \bar\nu')$ so that $\bar\nu' : M' \to B$ is an $[\frac{m}{2}]$-equivalence.*

*Proof.* In the first step we make $M$ connected. We can diminish the number of components of $M$ by one if we do surgery via an appropriate embedding $f : S^0 \times D^m \hookrightarrow M$, if $f(1, 0)$ and $f(-1, 0)$ are contained in different components of $M$ (note that in this situation surgery is the same as forming the connected sum).

Now, we assume $M$ to be connected and deal in the second step with $\pi_1$. We want to modify $\bar\nu : M \to B$ so that the induced map in $\pi_1$ is surjective.

For this, and the similar statement for higher homotopy groups, it is useful to note that surgery on a standard (unknotted) embedding $S^i \times D^{m-i} \hookrightarrow D^m \hookrightarrow M$ replaces $M$ by $M \# S^{i+1} \times S^{m-i-1}$. More precisely, consider the decomposition of $S^m = S^i \times D^{m-i} \cup D^{i+1} \times S^{m-i-1}$. Surgery on $S^i \times D^{m-i}$ yields $S^{i+1} \times S^{m-i-1}$ and replacing $M$ by $M \sharp S^m$ we obtain via surgery $M' = M \# S^{i+1} \times S^{m-i-1}$. We have freedom in extending the normal $B$-structure on $M$ to the trace of the surgery and this freedom can be used to achieve the fact that under the restriction of the normal $B$-structure on the trace to $M'$ an arbitrary element in the kernel of $\pi_{i+1}(B) \to \pi_{i+1}(BO)$ is in the image of $\bar{\nu}'_*$.

We can generalize this construction. For $\alpha : S^i \to O(m-i)$ twist the embedding of $S^i \times D^{m-i}$ by composition with the corresponding diffeomorphism on $S^i \times D^{m-i}$. Performing surgery replaces $M$ by $M \sharp X_\alpha$ where $X_\alpha$ is the sphere bundle of the vector bundle over $S^{i+1}$ classified by $\alpha$. If $\alpha \in \pi_i(BO)$ is in the image of $\pi_{i+1}(B) \to \pi_{i+1}(BO)$, the normal $B$-structure on $M$ extends to the trace of the surgery and now $\alpha$ is in the image of the map induced by the normal Gauss map from $M \sharp X_\alpha$ to $BO$.

We call such surgeries *connected sum surgeries*. Combining these two considerations and using the fact that $\pi_i(B)$ is finitely generated (over $\Lambda$ for $i > 1$), we obtain:

LEMMA 3. *For $i \leq m/2$, by a sequence of connected sum surgeries, $\bar{\nu}_* : \pi_i(M) \to \pi_i(B)$ is surjective without changing anything below dimension $i$.*

Let $\langle x_1, \dots, x_k | r_1, \dots, r_s \rangle$ be a presentation of $\pi_1(B)$. Applying the lemma above to $\pi_1$ we can replace $(M, \bar{\nu})$ by $M' = M \sharp X, \bar{\nu}'$ with $X$ a connected sum of $X'_\alpha$s as above, such that $\pi_1(M')$ has a presentation

$$\langle a_1, \dots, a_j, z_1, \dots, z_k | R_1, \dots, R_p \rangle,$$

where $\langle z_i \rangle = \pi_1(X)$, $\bar{\nu}'_* z_i = x_i$ and $\langle a_1, \dots, a_j | R_1, \dots, R_p \rangle$ is a presentation of $\pi_1(M)$ (note that by Morse theory [Mi1] the fundamental group of $M$ is finitely presentable if $M$ is compact, in particular $r_i$ is a word in $a_1, \dots, a_j$).

In this situation we write $\bar{\nu}'_*(a_i) = w_i(x_1, \dots, x_k)$, a word in $x_i$. Now consider the elements $a_i^{-1} w_i(z_1, \dots, z_k)$ in $\pi_1(M')$ and $r_i(z_1, \dots, z_k)$. Obviously these elements are in the kernel of $\bar{\nu}'_*$ and thus we can do surgery on them. The effect on $\pi_1(M')$ is to introduce these elements as additional relations. This follows from Lemma 2 since $m \geq 4$. Thus the map on $\pi_1$ becomes an isomorphism. By Lemma 3 we can assume that $\pi_2(M) \to \pi_2(B)$ is surjective.

Summarizing after these steps we can assume that $\bar{\nu} : M \to B$ with $M$ connected and $\bar{\nu}$ a 2-equivalence. We finish the proof by an inductive argument.

We assume inductively that for some $2 \leq r < [\frac{m}{2}]$, $\bar{\nu}$ is an $r$-equivalence. We first want to eliminate the kernel of $\bar{\nu}_* : \pi_r(M) \to \pi_r(B)$ by a sequence of surgeries. There is an exact sequence

$$\pi_{r+1}(B, M) \xrightarrow{d} \pi_r(M) @>\bar{\nu}_*>> \pi_r(B) \to 0$$

(here as in similar situations we replace $\bar{\nu} : M \to B$ by an embedding up to homotopy equivalence using the mapping cylinder, so that the relative homotopy groups make sense [Wh]).

By assumption, $B$ has a finite $(r + 1)$-skeleton so that $H_{r+1}(B, M; \Lambda) \cong \pi_{r+1}(B, M)$ is finitely generated. Surgery on a set of generators of image $d$ eliminates the kernel of $\bar\nu_*$ without changing the inductive assumptions (this follows from Lemmas 1 and 2). Finally, as for $r = 0$ and 1, we can do connected sum surgeries to show that $\pi_{r+1}(M) \to \pi_{r+1}(B)$ is surjective. $\qquad \square$

We call two compact manifolds $M_0$ and $M_1$ with the same boundary and normal $B$-structures, which agree on the boundary, *normally $B$-bordant relative to the boundary*, if the union of the two manifolds over the common boundary is zero bordant as a normal $B$-manifold. Here we have to equip $M_1$ with the negative orientation which is obtained by extending the given $B$-structure on $M_1$ to the cylinder $M_1 \times I$ and restricting it to the other boundary component.

Obviously the trace of a surgery is a normal $B$-bordism relative boundary. Thus, we can conclude from Proposition 3 the following:

COROLLARY 1. *Under the assumptions of Proposition 4, $(M, \bar\nu)$ is normally $B$-bordant relative to the boundary to $(M', \bar\nu')$ such that $\bar\nu' : M' \to B$ is an $[\frac{m}{2}]$-equivalence.*

The concept of normal 1-types and normal $B$-bordisms is useful for the investigations of a relevant differential geometric problem: Which manifolds admit a metric of positive scalar curvature? This relation was pointed out to me by Stephan Stolz. The key is the following result which is an easy consequence of the surgery theorem of Gromov-Lawson [G-L], respectively, Schoen-Yau [S-Y].

THEOREM 1 [G-L], [S-Y]. *Let $M$ be a compact manifold of dimension $n \geq 5$. Let $B$ be the normal 1-type of $M$ as described in Proposition 2. Then $M$ admits a metric of positive scalar curvature if and only if there is a normal $B$-manifold $N$ admitting a metric of positive scalar curvature, such that $M$ and $N$ agree in $\Omega_n(B)/\mathrm{Aut}(B)$.*

*Proof.* Let $(W, \bar\nu_W)$ be a normal $B$-bordism between $(M, \bar\nu)$ and $(M', \bar\nu')$. By Proposition 4 we can assume that $\bar\nu_W$ is a 3-equivalence, implying that $i : M \to W$ is a 2-equivalence. By Morse theory $M$ is obtained from $M'$ by a sequence of surgeries [Mi1]. If $i : M \to W$ is a 2-equivalence the proof of this theorem implies that one actually can pass from $M'$ to $M$ by a sequence of surgeries using embeddings of $S^r \times D^{m-r}$ with $r < m - 2$ [Mi2].

The surgery theorem of [G-L] or [S-Y] says that if one performs surgery on a sphere of codimension $\geq 3$ on a manifold with positive scalar curvature metric, then the resulting manifold admits such a metric. Thus the existence of a positive scalar curvature metric on $M'$ implies the existence on $M$. $\qquad \square$

Corollary 2.  *Let $M$ be a closed manifold of dimension $m \geq 5$ admitting a zero bordant normal 1-smoothing $\bar{\nu}$ in $\xi$ where $\xi$ is the normal 1-type of $M$ as described in Proposition 2. Then $M$ admits a metric with positive scalar curvature.*

*Proof.*  $(M, \bar{\nu})$ is $B$-bordant to the sphere $S^m$ with the normal $B$-structure induced from $D^{m+1}$. Since the standard metric on $S^m$ has positive sectional curvature (implying positive scalar curvature), the result follows from Theorem 1.  □

*Remark.* For $M$ simply connected of dimension $\geq 5$ the solution of the problem of existence of a positive scalar curvature metric follows if $M$ does not admit a spin structure ($w_2(M) \neq 0$). For, in this situation one can rather easily construct explicit generators of the oriented bordism group $\Omega_n$ (which in this situation is the bordism group of the normal 1-type) admitting metrics of positive scalar curvature. This was carried out in [G-L]. The spin case is much harder and was recently solved by Stephan Stolz [Sto1] showing that there is a single obstruction $\alpha(M)$ with values in $\mathbb{Z}$ for $\dim(M)$ divisible by 4, in $\mathbb{Z}/2$ for $\dim(M) \equiv 0, 1 \mod 8$ and 0 else. There is also substantial progress going on for nonsimply connected manifolds [Ro-St], [Sto2], [Ju].

## 4. Stable diffeomorphism classification

In this section we will prove Theorem C and a relative version for manifolds with boundary. We will do it by showing that a normal $B$-bordism $W$ between two normal $B$-smoothings of $2q$-dimensional manifolds $M_0$ and $M_1$ in a $(q-1)$-universal fibration $B$ can be replaced by an s-cobordism after a sequence of surgeries and a new operation, called subtraction of tori, which changes the boundary components by connected sum with $S^q \times S^q$. Then the s-cobordism theorem [Ke] in dimension $> 4$ and the stable s-cobordism theorem in dimension 4 [Q] imply that $M_0$ and $M_1$ are stably diffeomorphic.

We will also prove a relative version for manifolds with boundary. Let $M_0$ and $M_1$ be compact manifolds of dimension $2q$ with boundary and $f : \partial M_0 \to \partial M_1$ a diffeomorphism. This diffeomorphism is used to identify the boundaries. Suppose that these manifolds have the same normal $(q-1)$-type and admit normal $(q-1)$-smoothings compatible with $f$, i.e. are equal on the boundary after identifying the boundaries via $f$. We also assume that the normal $B$-manifold $M_0 \cup_f M_1$ is zero bordant via a normal $B$-bordism $W$. We begin with the description of subtraction of tori from $W$. Consider an embedded torus $S^q \times D^{q+1}$ in the interior of $W$. Join $\partial(S^q \times D^{q+1})$ with $M_0$ by an embedded thickened arc $I \times D^{2q}$ meeting $\partial(S^q \times D^{q+1})$ and $M_0$ transversely in $\{0\} \times D^{2q}$ and $\{1\} \times D^{2q}$ respectively. Remove $S^q \times \text{int}(D^{q+1})$ and $I \times \text{int}(D^{2q})$ from $W$ and straighten the resulting angles (compare [C-F, p. 9]). The resulting manifold $W'$ has boundary $M_0 \# S^q \times S^q \cup_f M_1$. We say that $W'$ is obtained from $W$ by *subtraction of a (solid) torus*. Of course, we can do the same with $M_1$ instead of $M_0$. One can generalize this process by admitting embeddings of arbitrary vector bundles over $S^q$ instead of the trivial bundle. Then one stabilizes by connected sum with the corresponding sphere bundle. Also this generalization is useful for some classification problems (compare [Kr1], [Kr2]).

We want to do this process with a bit more care controlling the $B$-structures. Up to homotopy classes of lifts $D^{q+1}$ has a unique normal structure in $B$ and we denote its restriction to $S^q$ by $\bar{\nu}_c$ (note that this "canonical" lift is not the constant map as $D^{q+1} \subseteq \mathbb{R}^{q+2}$ meets $\mathbb{R}^{q+1}$ transversely in $S^q$). Similarly, we can construct a normal structure on $S^q \times D^{q+1}$ and we denote its restriction to $S^q \times S^q$ again by $\bar{\nu}_c$. Now, we will show that, if $S^q \times \{0\}$ is zero homotopic in $B$, we can change the embedding of $S^q \times D^{q+1}$ into $W$, such that the restriction of the normal $B$-structure on $W$ to $M_0 \# S^q \times S^q$ is equal to $M_0 \# (S^q \times S^q, \bar{\nu}_c)$. For this, we note that the different normal $B$-structures on $S^q \times D^{q+1}$ are classified up to homotopy by the action of $\pi_q(F)$ on a fixed given normal $B$-structure, where $F$ is the fibre of $B \longrightarrow BO$. Since $S^q \times \{0\}$ is zero homotopic in $B$ we are only allowed to change the $B$-structure on $S^q \times \{0\}$ by elements in the image of $\pi_{q+1}(BO) \to \pi_q(F)$.

Now we consider a map $\alpha : S^q \longrightarrow O(q+1)$ and the twist diffeomorphism $f_\alpha : S^q \times D^{q+1} \longrightarrow S^q \times D^{q+1}$, $(x, y) \longrightarrow (x, \alpha(x) \cdot y)$.

The induced normal $B$-structure under this diffeomorphism on $S^q \times D^{q+1}$ is given by the action of the image of $\alpha$ under $\pi_q(O(q+1)) \longrightarrow \pi_q(O) \cong \pi_{q+1}(BO) \to \pi_q(F)$ on the given $B$-structure. Since $\pi_q(O(q+1)) \longrightarrow \pi_q(O)$ is surjective [Ste], this implies that we can always change a given embedding of $S^q \times D^{q+1}$ into $W$ by composing it with $f_\alpha$ for an appropriate $\alpha$ such that the induced normal $B$-structure on $\partial(S^q \times D^{q+1})$ is fibre homotopic to $\nu_c$. In the following we will always assume that embeddings $S^q \times D^{q+1}$ into $W$ are chosen such that the normal $B$-structure on $\partial(S^q \times D^{q+1})$ is $\nu_c$. Then we call this a *compatible subtraction of a torus*.

THEOREM 2. *Let $M_0$ and $M_1$ be compact connected $2q$-dimensional manifolds with normal $(q-1)$-smoothings in a fibration $B$. Let $f : \partial M_0 \to \partial M_1$ be a diffeomorphism compatible with the normal $(q-1)$-smoothings.*

*By a finite sequence of surgeries and compatible subtraction of tori, a normal $B$-zero bordism of $M_0 \cup_f M_1$ can be replaced by a relative s-cobordism between $M_0 \sharp r(S^q \times S^q)$ and $M_1 \sharp s(S^q \times S^q)$.*

COROLLARY 3. *Under the same conditions $f : \partial M_0 \to \partial M_1$ can be extended to a diffeomorphism $F : M_0 \sharp r(S^q \times S^q) \to M_1 \sharp s(S^q \times S^q)$. This diffeo-* *morphism commutes up to homotopy with the normal $(q-1)$-smoothings in $B$ given by the normal $(q-1)$-smoothing on $M_i$ and $\nu_c$ on $S^q \times S^q$.*

If the manifolds have the same Euler characteristic, then $r = s$. If the boundary is empty, this is Theorem C from the introduction.

*Proof.* In the following we will frequently make use of homology and cohomology with twisted coefficients. In particular, we consider as coefficients the group ring $\Lambda = \mathbb{Z}[\pi_1(B)]$. Here we assume that the space whose homology we are looking at is equipped with a map to $B$ under which we pull back the coefficients. In particular, if the map induces an isomorphism on $\pi_1$, the homology with coefficients in $\Lambda$ is the ordinary homology of the universal covering considered as a module over $\pi_1$ via covering translations. Note that the corresponding statement for cohomology is only true for finite $\pi_1$; for infinite groups it is ordinary cohomology with compact support. References for homology with twisted coefficients are [Wh], [W1].

$W$ is a relative $s$-cobordism if and only if

i) $\pi_1(M_i) \longrightarrow \pi_1(W)$ are isomorphisms for $i = 0, 1$.

ii) $H_k(W, M_i; \Lambda) = \{0\}$ for $i = 0, 1$ and $k \leq q$.

iii) The Whitehead torsion $\tau(W, M_i)$ vanishes for $i = 0, 1$ [Mi3].

By Proposition 4 we can assume that $\bar{\nu} : W \longrightarrow B$ is a $q$-equivalence. Since also the normal $(q-1)$-smoothings $\bar{\nu}_i : M_i \to B$ are $q$-equivalences, this implies i) and that ii) holds for $k < q$ . To kill $H_q(W, M_i; \Lambda)$ by a sequence of compatible subtractions of tori, we consider the diagram of exact sequences

$$
\begin{array}{ccccccccc}
 & & H_{q+1}(B, W; \Lambda) & & & & & & \\
 & & \downarrow & & & & & & \\
H_q(M_i; \Lambda) & \to & H_q(W; \Lambda) & & \to & H_q(W, M_i; \Lambda) & \to & 0. \\
 & & \downarrow & & & & & & \\
 & & H_q(B; \Lambda) & & & & & & \\
 & & \downarrow & & & & & & \\
 & & 0 & & & & & &
\end{array}
$$

As $H_q(M_i; \Lambda) \to H_q(B; \Lambda)$ is surjective, the same follows for $H_{q+1}(B, W; \Lambda) \to H_q(W, M_i; \Lambda)$. Since $W$ and $M_i$ are compact, $H_q(W, M_i; \Lambda)$ is a finitely generated $\Lambda$-module. As $\bar{\nu} : W \to B$ is a $q$-equivalence, the Hurewicz theorem implies

$$\pi_{q+1}(B, W) @ > \cong > > H_{q+1}(B, W; \Lambda).$$

Thus there exists a finite set of elements of $\pi_{q+1}(B, W)$ mapping to generators of $H_q(W, M_0; \Lambda)$. The images of them in $\pi_q(W)$ can be represented by disjointly embedded spheres with trivial normal bundle $(S^q \times D^{q+1})_i$ in the interior of $W$ (the normal bundle is stably trivial since these elements map to zero in $\pi_q(B)$ under $\bar{\nu}$ and we are in the stable range, i.e. stably trivial bundles are trivial). As described above, join each of these embedded $(S^q \times D^{q+1})_i$ with $M_0$ to obtain $[(S^q \times D^{q+1}) \cup I \times D^{2q}]_i$ and subtract these tori to obtain $(W', \bar{\nu}')$ replacing $(M_0, \bar{\nu}_0)$ by $(M_0, \bar{\nu}_0) \# r(S^q \times S^q, \bar{\nu}_c)$. This leaves $(M_1, \bar{\nu}_1)$ unchanged. By general position $\bar{\nu}'$ is again a $q$-equivalence. The fact that the classes $(S^q \times \{0\})_i$ generate $H_q(W, M_0; \Lambda)$ and the long exact homology sequence of a triple implies

$$H_q(W, M_0 \cup_i [(S^q \times D^{q+1}) \cup I \times D^{2q}]_i; \Lambda) = \{0\}.$$

But the latter group is by excision isomorphic to $H_q(W', M_0 \# r(S^q \times S^q); \Lambda)$.

Thus we have killed $H_q(W', M_0 \# r(S^q \times S^q); \Lambda)$ and so the pair $(W', M_0 \# r(S^q \times S^q))$ is $q$-connected. What about $H_q(W', M_1; \Lambda)$? This $\Lambda$-module is stably free; i.e., the direct sum with an appropriate finitely generated free $\Lambda$-module is free: It was shown in [W1, Lemma 2.3] that this is true if $H_k(W', M_1; A)$ is trivial for $k \neq q$ and every $\Lambda$-module $A$. This holds for $k < q$ by assumption and for $k > q$ we have by Poincaré-Lefschetz duality (compare [W1, Th. 2.1])

$$H_k(W', M_1; A) \cong H^{2q+1-k}(W', M_0 \# r(S^q \times S^q); A) = \{0\}$$

since $(W', M_0 \# r(S^q \times S^q))$ is $q$-connected.

We want to achieve that $H_q(W', M_1; \Lambda)$ is actually free. For this we consider a finite set of disjoint embeddings of $S^q \times D^{q+1}$ sitting unknotted and unlinked in a ball $D^{2q+1} \subset \text{int}(W)$. If we join them with $M_0$ and subtract these tori, we replace $H_q(W', M_1; \Lambda)$ by the direct sum with a free module of rank the number of embedded $S^q \times D^{q+1}$.

Thus we can assume that $H_q(W', M_1; \Lambda)$ is free. If one distinguishes a basis of $H_q(W', M_1; \Lambda)$, the Whitehead torsion $\tau(W', M_1)$ is defined [Mi3, p. 378]. Since every element of $\text{Wh}(\pi_1)$ is represented by a matrix of finite rank, the definition of this torsion implies that (after further stabilization) we can choose the basis such that the torsion vanishes. Following Wall, such a basis is called *preferred* [W1, §2].

Given such a preferred basis, we have shown above that we can represent it by $s$ disjoint embeddings $(S^q \times D^{q+1})_j$ in the interior of $W'$. Join them with $M_1$ and subtract these tori to obtain $W''$ with

$$\partial W'' = M_0 \# r(S^q \times S^q) \cup_f M_1 \# s(S^q \times S^q).$$

We claim that $W''$ is a relative $s$-cobordism. For this we have to check

$$H_q(W'', M_0 \# r(S^q \times S^q); \Lambda) = 0 = H_q(W'', M_1 \# s(S^q \times S^q); \Lambda)$$

and the torsion of $(W'', M_i)$ is zero.

As above, Poincaré-Lefschetz duality implies that the left homology group vanishes automatically if $H_k(W'', M_1 \# s(S^q \times S^q); \Lambda) = \{0\}$ for all $k$. Furthermore we know that these homology groups are trivial for $k \neq q$ or $q + 1$. By excision the group is isomorphic to $H_k(W', M_1 \cup U; \Lambda)$ where $U = \bigcup_j (S^q \times D^{q+1})_j \cup I \times D^{2q}$. We consider the homology sequence of a triple (with $\Lambda$-coefficients)

$$0 \to H_{q+1}(W', M_1 \cup U) \to H_q(M_1 \cup U, M_1)$$
$$\to H_q(W', M_1) \to H_q(W', M_1 \cup U) \to 0.$$

All other homology groups in this long exact sequence vanish. $H_q(M_1 \cup U, M_1; \Lambda)$ is free with its basis given by the $(S^q \times \{0\})_j$. This is a preferred basis (the Whitehead torsion of $(M_1 \cup U, M_1)$ with respect to it vanishes). The image of these basis elements forms our preferred basis of $H_q(W', M_1)$.

All this implies that $H_k(W', M_1 \cup U; \Lambda) = \{0\}$ for all $k$ and that the Whitehead torsion of the based acyclic complex given by the exact sequence above is trivial. By the additivity formula of the Whitehead torsion [Mi3, Th. 3.2], this implies:

$$\tau(W'', M_1 \# s(S^q \times S^q)) = \tau(W', M_1 \cup U) = \tau(W', M_1) - \tau(M_1 \cup U, U) = 0.$$

By the duality theorem [Mi3], this also implies $\tau(W'', M_0 \# r(S^q \times S^q)) = 0$, finishing the proof of Theorem 2.  □

As Peter Teichner pointed out to me, one can get an easier proof of Corollary 3 roughly as follows (a similar argument was used in the proof of [Fr1, Th. 3]). By Proposition 4 we can assume that there is a bordism relative to the boundary $W$ between $M_0$ and $M_1$ such that the map $W \to B$ is a $q$-equivalence. Then the proof of the $s$-cobordism theorem implies that $W$ has a handle decomposition consisting only of handles of index $q$ and $q + 1$. The middle level of this bordism is obtained from both $M_0$ and $M_1$ as the result of surgeries on null-homotopic embeddings $S^{q-1} \times D^{q+1}$ in $M_i$ which replace $M_i$ by connected sum with copies of $S^q \times S^q$. Thus $M_0$ and $M_1$ are stably diffeomorphic relative to the boundary.

## 5. The main theorem for even-dimensional $W$

We begin with the definition of the obstruction monoid $l_{2q}(\pi, w)$. Here $\pi$ is a group and $w : \pi \to \mathbb{Z}/2$ a homomorphism, which in the geometric context is the first Stiefel-Whitney class. Denote the integral group ring as before by $\Lambda = \mathbb{Z}[\pi]$. Let $- : \Lambda \to \Lambda$ be the involution sending $g \in \pi$ to $\bar{g} = w(g) \cdot g^{-1}$. We will work with left $\Lambda$-modules but note that we can via $-$ pass to right $\Lambda$-modules whenever we like. In particular, the dual $\Lambda$-module $V^\star$ of a left $\Lambda$-module $V$ is naturally a right $\Lambda$-module but we consider it as left $\Lambda$-module.

For $\varepsilon = \pm 1$, an *$\varepsilon$-quadratic form* over $\Lambda$ is given by a left $\Lambda$-module $V$ together with an $\varepsilon$-hermitian form $\lambda : V \times V \to \Lambda$ and a quadratic refinement $\mu : V \to \Lambda/_{\langle x - \varepsilon \bar x \rangle}$. This means that $\lambda$ and $\mu$ have to fulfill the following properties:

i) For fixed $v \in V$ the map $V \to \Lambda$ mapping $w \mapsto \lambda(w, v)$ is a $\Lambda$-homomorphism.

ii) $\lambda(v, w) = \varepsilon \bar\lambda(w, v)$.

iii) $\lambda(v, v) = \mu(v) + \varepsilon \bar\mu(v) \in \Lambda$.

iv) $\mu(v + w) = \mu(v) + \mu(w) + \lambda(v, w) \in \Lambda/_{\langle x - \varepsilon \bar x \rangle}$.

v) $\mu(x \cdot v) = x \cdot \mu(v) \cdot \bar x$.

Note that iii) is an equation in $\Lambda$ since $x + \varepsilon \bar x$ is a well-defined element in $\Lambda$. An important special case is the *$\varepsilon$-hyperbolic form* $H^r_\varepsilon := H_\varepsilon \perp \cdots \perp H_\varepsilon$, $r$ summands, where $H_\varepsilon$ is the form on $\Lambda \oplus \Lambda$ with standard basis $e$ and $f$ and $\lambda(e, f) = 1$, $\lambda(f, e) = \varepsilon$, $\lambda(e, e) = \lambda(f, f) = 0$ and $\mu(e) = \mu(f) = 0$.

A $\Lambda$-module $V$ is called *based*, if it is finitely generated and equipped with an equivalence class of bases, where two bases are equivalent if the matrix of base changes vanishes in $\mathrm{Wh}(\pi)$. An isomorphism between based $\Lambda$-modules is called a *simple isomorphism* if the matrix of the isomorphism with respect to the given bases vanishes in $\mathrm{Wh}(\pi)$.

The objects in $l_{2q}(\pi, w)$ are represented by triples $(V_0 \xleftarrow{f} V \xrightarrow{g} V_1, \lambda, \mu)$ fulfilling:

i) $V$ is a finitely generated left $\Lambda$-module and $V_0$ and $V_1$ are based.

ii) $\lambda : V_0 \to V_1^*$ is an isomorphism and should induce an $\varepsilon$-hermitian form $\lambda_V$ on $V$ (i.e. for $x$ and $y$ in $V$ we have $\lambda_V(x, y) := \lambda(f(x), g(y)) = (-1)^q \overline{\lambda(f(y), g(x))} \in \Lambda$) and $\mu : V \to \Lambda/_{\langle x - \varepsilon \bar x \rangle}$ is a quadratic refinement of this form. Here $\varepsilon := (-1)^q$.

The second condition can be reformulated as follows: The adjoint of $\lambda$ is an isomorphism $\lambda : V_0 \to V_1^*$ and $\lambda_V = g^* \lambda f = (-1)^q f^* \lambda^* g : V \to V^*$, and $\mu : V \to \Lambda/_{\langle x - \varepsilon \bar x \rangle}$ is a quadratic refinement of $\lambda_V$.

The orthogonal sum defines a monoid structure on these objects. Particular objects of this type are obtained from an $\varepsilon$-quadratic form $(V, \lambda, \mu)$ with $V$ based by setting $V_0 = V_1 = V$ and $f = g = \mathrm{Id}$. In particular, the $\varepsilon$-hyperbolic forms are of this type, where we base $H_\varepsilon$ via the canonical basis $e$ and $f$. In the definition of the quadratic refinement we can, for $\varepsilon = -1$, replace $\mu$ by $\tilde\mu$, which takes values in the quotient $\Lambda/\langle x - \varepsilon \bar x, 1 \rangle$. Everything above makes sense with this modification.

*Definition.* The monoid $l_{2q}(\pi, w)$ is given by equivalence classes of triples

$$(V_0 \xleftarrow{f} V \xrightarrow{g} V_1, \lambda, \mu)$$

as above, where two such triples are equivalent if they are simple isometric after adding $\varepsilon$-hyperbolic forms (of perhaps different rank) to them. If, for $q$ odd, we replace $\mu$ by $\tilde{\mu}$ which takes values in $\Lambda/\langle x - \varepsilon \bar{x}, 1 \rangle$, we call the corresponding monoid $l^{\sim}_{2q}(\pi, w)$.

Of course, $l_{2q}(\pi, w)$ depends only on $q$ modulo 2. The notion comes from the geometric context. There is a submonoid coming from all $\varepsilon$-quadratic forms $(V, \lambda, \mu)$ with $V$ based. This is actually an *abelian group* $L_{2q}(\pi, w)$ which essentially is the ordinary Wall group $L^s_{2q}(\pi, w)$. More precisely, the Wall group is given by those $\varepsilon$-quadratic forms $(V, \lambda, \mu)$ with $V$ based, where the adjoint of $\lambda$ is a simple isomorphism. Thus we have a homomorphism from the subgroup of $l_{2q}(\pi, w)$ consisting of $\varepsilon$-quadratic forms $(V, \lambda, \mu)$ with $V$ based to $\mathrm{Wh}(\pi)$ mapping to the adjoint of $\lambda$ with kernel the *Wall group* $L^s_{2q}(\pi, w)$.

Now we recall Wall's definition of a quadratic form on even-dimensional manifolds. Let $W$ be a $2q$-dimensional compact manifold. We equip $W$ with a base point and orient it at this base point. Wall [W1, p. 44ff] defines a $(-1)^q$-hermitian form (i.e. the conditions i) and ii) above are fulfilled) on the group of regular homotopy classes of immersions of $q$-spheres in $W$, denoted $\mathrm{imm}_k(W)$. Roughly speaking, the hermitian form is given by transversal double point intersections which, along the two branches, are joined with the base points, so that the form has values in $\Lambda$. Similarly, he assigns to each immersion $x$ an element $\mu(x) \in \Lambda/\langle x - (-1)^q \bar{x} \rangle$ which is given analogously via self-intersection points. We recall this fundamental result.

PROPOSITION 5 ([W1, Th. 5.2 and p. 52]). *Let $W$ be a compact $2q$-dimensional manifold with base point and local orientation at the base point.*

*Intersections define a $(-1)^q$-hermitian form $\lambda : \mathrm{imm}_k(W) \times \mathrm{imm}_k(W) \to \Lambda$ and self-intersections define a function $\mu : \mathrm{imm}_k(W) \to \Lambda/\langle x - (-1)^q \bar{x} \rangle$. If $q \geq 3$, an immersion is regularly homotopic to an embedding if and only if $\mu$ vanishes on it.*

*If $W$ is part of the boundary of some manifold $X$ with same fundamental group and there are two families of disjoint immersions of spheres $S^k$ in $W$ each of which extends to an immersion of a disk with holes, then the sum of the self-intersection of all immersed spheres within each family and the intersection number between the two families vanish (slogan: intersection and self-intersection numbers vanish for elements in the kernel of $\mathrm{imm}_k(W) \to \mathrm{imm}_k(X)$).*

The function $\mu$ is not a quadratic refinement of $\lambda$ but very closely related to a quadratic refinement. It fulfills for $v, w \in \mathrm{imm}_q(W)$:

iii') $\lambda(v,v) = \mu(v) + \varepsilon\bar\mu(v) + e(v) \in \Lambda$, $e$ the normal Euler number.

iv') $\mu(v+w) = \mu(v) + \mu(w) + \lambda(v,w) \in \Lambda/\langle x - \varepsilon\bar x\rangle$.

v') $\mu(g \cdot v) = g \cdot \mu(v) \cdot \bar g$ for $g \in \pi_1(W)$.

In [W1, Th. 5.2] the last formula was stated for arbitrary $a \in \Lambda$ instead of $g \in \pi_1(W)$, which is not correct. But the proof shows v') and is correct. This function $\mu$ induces a quadratic refinement on those homotopy classes which have stably trivial normal bundle. To explain this, we first study when an immersion is homotopic to an embedding with trivial normal bundle. A necessary condition is that the stable normal bundle of the immersion is trivial or equivalently, that the homotopy class $\alpha$ represented by the immersion is contained in $K\pi_q(W) := Ker\nu_* : \pi_q(W) \to \pi_q(BO)$. For $q$ even, a stably trivial $q$-dimensional bundle over $S^q$ is trivial if and only if the Euler class vanishes and the group of these bundles is generated by the tangent bundle of $S^q$. Since the Euler class is controlled by the self-intersection number, one can for $q$ even add in the third sentence of the proposition that the embedding has trivial normal bundle, if and only if the stable normal bundle is trivial. For $q \neq 1, 3, 7$ odd there are precisely two stably trivial bundles of dimension $q$ over $S^q$, the trivial bundle and the tangent bundle of $S^q$ [K-M, p. 534]. For $q = 1, 3, 7$, $q$-dimensional stably trivial vector bundles over $S^q$ are trivial.

The map $\mathrm{imm}_q(W) \to \pi_q(W)$ is surjective and its kernel is in the image of a homomorphism from $\mathbb{Z}$ or $\mathbb{Z}/2$ to $\mathrm{imm}_q(W)$ if $q$ is even or odd [W1, p. 44]. Assume that $W$ decomposes as $W'\sharp S^q \times S^q$. Then there are two immersions representing the diagonal of $S^q \times S^q$, the diagonal embedding and the connected sum of the two factors. The normal bundle of the first is the tangent bundle and the normal bundle of the second is trivial. The normal bundle of the difference of these two immersions is the tangent bundle of $S^q$. Thus, after adding an appropriate multiple of the difference of these two immersions, we can assume that an arbitrary immersion representing a homotopy class in $K\pi_q(W)$ is represented by an immersion with trivial normal bundle and, for $q \neq 1, 3, 7$, this immersion is unique. This gives for $q \neq 1, 3, 7$ a splitting of the restriction of $\mathrm{imm}_q(W) \to \pi_q(W)$ to the kernel of the map to $\pi_q(BO)$. We can now define $\mu$ on $K\pi_q(W)$. By the formulas iii'), iv') and v') this is a quadratic refinement of $\lambda$ and $\alpha \in K\pi_q(W)$ is representable by a compatible embedding with trivial normal bundle if and only if $\mu(\alpha) = 0$. If $q = 3, 7$, there is no unique immersion in $K\pi_q(W)$ with trivial normal bundle. The difference of the diagonal in $S^q \times S^q$ and the sum of the two factors is the nontrivial element in the kernel of $\mathrm{imm}_q(W) \to \pi_q(W)$. It is nontrivial since $\mu$ is 1 on it. Thus the self-intersection is not well-defined in this case. If there is $\beta \in \pi_{r+1}(B)$ with $\beta^*\xi^*w_{r+1} \neq 0$, $w_{r+1} \in H^*(BO; \mathbb{Z}/2)$ the Stiefel Whitney class, we replace $\mu$ by $\tilde\mu$. Then by Lemma 2 we have again that $\alpha \in K\pi_q(W)$ is representable by a compatible embedding with trivial normal bundle if and only if $\tilde\mu(\alpha) = 0$.

We summarize these considerations as:

PROPOSITION 6. *For $q \neq 1, 3, 7$ the intersections and self-intersections define a quadratic form $(\lambda, \mu)$ on $K\pi_q(W)$. The statements of the previous proposition hold for this quadratic form. In particular, $\alpha \in K\pi_q(W)$ is representable by a compatible embedding with trivial normal bundle if and only if $\mu(\alpha) = 0$. The same holds for $q = 3, 7$ and $\langle w_{q+1}(B), \pi_{q+1}(B) \rangle \neq 0$ if $\mu$ is replaced by $\tilde{\mu}$.*

Now, let $M_0$ and $M_1$ be $(2q - 1)$-dimensional manifolds and $f : \partial M_0 \to \partial M_1$ be a diffeomorphism. Suppose that there are normal $(q - 2)$-smoothings in a fibration $B$ over $BO$ compatible with $f$. Let $W$ together with a normal $B$-structure $\bar{\nu}$ be a $B$-zero bordism of $M_0 \cup_f M_1$. Then we can by Proposition 4 assume that $\bar{\nu}$ is a $q$-equivalence, $W = W' \natural S^q \times S^q$ and by Hurewicz's theorem $H_{q+1}(B, W; \Lambda) \cong \pi_{q+1}(B, W)$. Thus the image under the boundary map $d : \pi_{q+1}(B, W) \to \pi_q(W)$ is a finitely generated $\Lambda$-module. This image is contained in $K\pi_q(W)$ and thus the quadratic form $(\lambda, \mu)$ is defined for $q \neq 1, 3, 7$. For $q = 1, 3, 7$ and $\langle w_{q+1}(B), \pi_{q+1}(B) \rangle \neq 0$ we replace $\mu$ by $\tilde{\mu}$.

For $q = 1, 3, 7$ and $\langle w_{q+1}(B), \pi_{q+1}(B) \rangle = 0$, we will again define a quadratic refinement $\mu$, but this can only be defined on $\mathrm{im}(d : \pi_{q+1}(B, W) \to \pi_q(W)) = \mathrm{Ker}(\pi_q(W) \to \pi_q(B))$, as follows. Let $\alpha \in \mathrm{Ker}(\pi_q(W) \to \pi_q(B))$ be represented by an immersion of $S^q$ into $W$. This has stably trivial and thus trivial normal bundle, since $q = 1, 3, 7$. Thus we can extend it to an immersion of $S^q \times D^q$ to $W$. The extension can be modified by twisting with elements of $\pi_q(0(q))$ as discussed in the proof of Lemma 2. The immersion of $S^q \times D^q$ to $W$ induces, together with the framing of $\nu(W)|_{\mathrm{im}(S^q)}$ given by the normal $B$-structure, a stable framing on $\nu(S^q)$.

We say that the immersion of $S^q \times D^q$ to $W$ is *good* if we can choose an extension to $S^q \times D^q$ in $W$ such that this stable framing is the standard one coming from the embedding $S^q \subset D^{q+1}$. The good immersions form a subgroup of $\mathrm{Ker}(\mathrm{imm}_q(W) \to \pi_q(W) \to \pi_q(B))$. It turns out that this subgroup maps isomorphically to $\mathrm{Ker}(\pi_q(W) \to \pi_q(B))$. The reason is that, as discussed above, $\mathrm{Ker}(\mathrm{imm}_q(B) \to \pi_q(B)) = \mathbb{Z}/2$ is generated by the difference of the diagonal in $S^q \times S^q$ and the immersion given by the connected sum of the two factors. The latter immersion is good while the diagonal is not as we will show below. Since $\pi_q(0)/\pi_q(0(q)) \cong \mathbb{Z}/2$ (see proof of Lemma 2), we can for every element $x$ in $\mathrm{Ker}(\pi_q(W) \to \pi_q(B))$ find a unique good immersion representing it. By this unique good immersion we define $\mu(x)$ using Proposition 5. Again $\mu(x) = 0$ if and only if $x$ is representable by a compatible embedding of $S^q \times D^q \hookrightarrow W$.

To show that the diagonal $\Delta$ in $S^q \times S^q$ (for $q = 1, 3, 7$) admits no compatible framing, we note that the normal bundle of $\Delta$ in $S^q \times S^q$ is the tangent bundle of $S^q$ and no unstable framing of the tangent bundle extends to $D^{q+1}$. But the restriction of the framing on $S^q \times S^q$ to $\Delta$ is the sum of two equal framings. Since the obstruction for a compatible framing is in $\pi_q(0)/_{\pi_q(0(q+1))} \cong \mathbb{Z}/2$ this implies that there is no compatible framing on the diagonal.

We summarize these considerations as:

PROPOSITION 7. *For $q = 3, 7$ and $\langle w_{q+1}(B), \pi_{q+1}(B)\rangle = 0$, there is a quadratic refinement $\mu$ of $\lambda$ defined on $\mathrm{im}(d : \pi_{q+1}(B, W) \to \pi_q(W)) = \mathrm{Ker}(\pi_q(W) \to \pi_q(B))$, and $\alpha \in K\pi_q(W)$ is representable by a compatible embedding with trivial normal bundle if and only if $\mu(\alpha) = 0$.*

Now we note that Poincaré duality gives a unimodular pairing

$$\lambda : H_q(W, M_0; \Lambda) \to H_q(W, M_1; \Lambda)^\star.$$

As in Section 4, we conclude from Poincaré duality that $H_q(W, M; \Lambda)$ are stably free. After stabilization of $W$ with $S^q \times S^q$'s we assume that these $\Lambda$-modules are free and equipped with a basis, so that the Whitehead torsion of $(W, M_i)$ vanishes. By Hurewicz's isomorphism we identify $\pi_q(W, M_i)$ with $H_q(W, M_i; \Lambda)$. The inclusions define homomorphisms $f$ and $g$:

$$\mathrm{im}(d : \pi_{q+1}(B, W) \to \pi_q(W)) \to H_q(W, M_i; \Lambda))$$

and the geometric interpretation of Poincaré duality implies that the intersection pairing $\lambda$ induces on $\mathrm{im}(d : \pi_{q+1}(B, W) \to \pi_q(W))$ the hermitian form as considered above.

*Definition.* For $q \neq 1, 3, 7$, or $q = 1, 3, 7$ and $\langle w_{q+1}(B), \pi_{q+1}(B)\rangle = 0$, we define, for $(W, \bar\nu)$ with the properties above, $\theta(W, \bar\nu) \in l_{2q}(\pi_1(B), w_1(B))$ by

$$\theta(W, \bar\nu) := (H_q(W, M_0; \Lambda) \xleftarrow{f} \mathrm{im}(d : \pi_{q+1}(B, W)$$
$$\to \pi_q(W)) \xrightarrow{g} H_q(W, M_1; \Lambda), \lambda, \mu),$$

where $f$ and $g$ are induced by inclusions.

For $q = 1, 3, 7$ and $\langle w_{q+1}(B), \pi_{q+1}(B)\rangle \neq 0$ we replace $\mu$ by $\tilde\mu$ in this definition and obtain $\theta(W, \bar\nu) \in l_{\widetilde{2q}}(\pi_1(B), w_1(B))$.

Here we identify $\pi_1(W)$ with $\pi_1(B)$ and $w_1(B)$ with $w_1(W) : \pi_1(W) \cong \pi_1(B) \to \mathbb{Z}/2$. Let $(W', \bar\nu)$ be another normal $B$-manifold with the same properties, which is normally $B$-bordant to $(W, \bar\nu)$ relative to the boundary. Then by Corollary 3, $W$ and $W'$ are stably diffeomorphic relative to the boundary and thus

$$\theta(W, \bar\nu) = \theta(W', \bar\nu') \in l_{2q}(\pi_1(B), w_1(B)) \ (\mathrm{resp.} l_{\widetilde{2q}}(\pi_1(B), w_1(B))).$$

Thus we can define $\theta(W, \bar\nu)$ for arbitrary $B$ zero-bordisms $(W, \bar\nu)$ of $M_0 \cup_f -M_1$ by replacing $W$ by a bordism with the properties above and then applying the definition above. This invariant $\theta(W, \bar{\nu})$ depends only on the bordism class relative to the boundary of $(W, \bar{\nu})$.

*Remark.* If the smoothings on $M_i$ are normal $k$-smoothings for some $k \geq q$, the invariant is contained in the subgroup $L_{2q}(\pi_1(B), w_1(B))$ and if, in addition, $B$ is a finite simple Poincaré complex and the smoothings are simple homotopy equivalences the obstruction sits in $L_{2q}^s(\pi_1(B), w_1(B))$. The reader might wonder, if the condition finitely generated free in the definition of elements in $l_{2q}(\pi, w)$ is necessary in the geometric context. Perhaps this is automatically the case for geometrically realized elements. Unfortunatley this is not the case, even if the bordism $W$ is already an s-cobordism. Consider for example a 3-dimensional lens space $L$ with nontrivial fundamental group $\pi$ and consider $W := L \times S^2 \times I$. Then the normal 2-type is $K(\pi, 1) \times \mathbb{CP}^\infty \times B\,\mathrm{Spin}$ and the obstruction is $(0 \leftarrow \pi_3(W) = \mathbb{Z} \times \mathbb{Z} \to 0, 0, 0)$, where $\mathbb{Z} \oplus \mathbb{Z}$ is the trivial $\pi$-module.

*Definition.* We call an element $\theta \in l_{2q}(\pi, w)$ (resp. $\theta \in \widetilde{l_{2q}}(\pi, w)$) *elementary* if and only if there is a representative $(V_0 \xleftarrow{f} V \xrightarrow{g} V_1, \lambda, \mu)$ and a based submodule $U \subset V$ such that

i) $U \subset U^\perp$, $\mu_{|U} = 0$ (or $\tilde{\mu}_{|U} = 0$, if $\theta \in \widetilde{l_{2q}}(\pi, w)$);

ii) $U$ maps injectively into $V_i$ and the image is a direct summand denoted $U_i$ whose basis is the image of the basis on $U$;

iii) $\lambda$ induces a simple isomorphism $U_0 \to (V_1/U_1)^*$.

Now, we can formulate the main theorem for $W$ even-dimensional.

THEOREM 3. *Let $M_0$ and $M_1$ be connected $(2q-1)$-dimensional manifolds, $q \geq 3$, and let $f : \partial M_0 \to \partial M_1$ be a diffeomorphism. Suppose that there are normal $(q-2)$-smoothings in a fibration $B$ over $BO$ compatible with $f$. Let $W$ together with a normal $B$-structure $\bar{\nu}$ be a $B$-zero bordism of $M_0 \cup_f M_1$. Then $(W, \bar{\nu})$ is $B$-bordant relative to the boundary of a relative s-cobordism if and only if $\theta(W, \bar{\nu}) \in l_{2q}(\pi_1(B), w_1(B))$ (resp. $\in \widetilde{l_{2q}}(\pi_1(B), w_1(B))$, if $q = 3, 7$ and $\langle w_{q+1}(B), \pi_{q+1}(B) \rangle \neq 0$) is elementary.*

*Proof.* By Proposition 4 we can assume that $\bar{\nu}$ is a $q$-equivalence; i.e. $(W, \bar{\nu})$ is a normal $(q-1)$-smoothing. Then $\theta(W, \bar{\nu})$ is defined as in the definition above. Since we can realize stabilization of $\theta(W, \bar{\nu})$ by hyperbolic forms geometrically via connected sum with $(S^q \times S^q)$'s, we can assume that for $\theta(W, \bar{\nu}) = (H_q(W, M_0; \Lambda) \longleftarrow \mathrm{im}(d : \pi_{q+1}(B, W) \to \pi_q(W)) \longrightarrow H_q(W, M_1; \Lambda), \lambda, \mu)$ there exists a based submodule $U \subset \mathrm{im}(d : \pi_{q+1}(B, W) \to \pi_q(W))$ with the properties i)–iii) in the definition of "elementary" above. Choose $x_1, \dots, x_k \in U$ representing the basis of $U$ implying $\lambda(x_i, x_j) = 0$ and $\mu(x_i) = 0$ (resp. $\tilde{\mu}(x_i) = 0$). By Propositions 6 and 7 one can find embeddings $(S^q \times D^q)_i$ representing $x_i$. Using the Whitney trick we can in addition assume that the embeddings are disjoint. Furthermore, we can assume that the embeddings are compatible with the $B$-structure so that we can do surgery with them (Lemma 2 and Propositions 6 and 7). We claim that the resulting $B$-manifold $(W', \overline{\nu}_{W'})$ is an $s$-cobordism.

It is clear that $\pi_1(M_i) \longrightarrow \pi_1(W')$ is an isomorphism. To compute $H_*(W', M_1; \Lambda)$ we consider the following exact sequences with $\Lambda$-coefficients. Write $X = \cup (S^q \times D^q)_i$

$$
\begin{array}{ccccc}
H_{q+1}(W', W - \overset{\circ}{X}) & & & & \\
\downarrow d & & & & \\
0 \to H_q(W - \overset{\circ}{X}, M_1) \to H_q(W, M_1) \overset{j}{\to} H_q(W, W - \overset{\circ}{X}) \to H_{q-1}(W - \overset{\circ}{X}, M_1) \to 0. \\
\downarrow & & & \downarrow \cong \\
H_q(W', M_1) & & & H_{q-1}(W', M_1) \\
\downarrow & & & \\
0 & & &
\end{array}
$$

By excision, $H_r(W, W - \overset{\circ}{X}; \Lambda) \cong H_r(X, \partial X; \Lambda)$. Thus $H_r(W, W - \overset{\circ}{X}; \Lambda)$ is trivial except for $r = q$ where it is $\Lambda^k$ with basis $[(\{*\} \times D^q, \{*\} \times S^{q-1})]_i$ or $r = 2q$ where it is again $\Lambda^k$ with basis $[(S^q \times D^q)_i, \partial]$. These bases together represent a preferred basis of $H_*(W, W - \overset{\circ}{X})$. We note that the map $U_0 \to H_q(W, W - \overset{\circ}{X})$ mapping $x_i \to [(\{*\} \times D^q, \{*\} \times S^{q-1})]_i$ is an isomorphism. Similarly, $H_*(W'W - \overset{\circ}{X}; \Lambda)$ has a preferred basis represented by $[(D^{q+1} \times \{*\}, \partial)]_i$ in dimension $q + 1$ and by $[(D^{q+1} \times S^{q-1})_i, \partial]$ in dimension $2q$.

With respect to this basis the homomorphism

$$H_q(W, M_1; \Lambda) \overset{j}{\longrightarrow} H_q(W, W - \overset{\circ}{X}; \Lambda)$$

is given by $x \longrightarrow (\lambda(x, x_1), \dots, \lambda(x, x_k))$. If we denote the image of $U$ in $H_q(W, M_i; \Lambda)$ by $U_i$, the definition of elementary implies that $H_q(W, M_1; \Lambda)$ splits as $U_1 \oplus H_q(W, M_1; \Lambda)/U_1$ and that the restriction of $j$ to $H_q(W, M_1; \Lambda)/U_1$ is a simple isomorphism to $U_0 \cong H_q(W, W - \overset{\circ}{X}; \Lambda)$. Thus $H_{q-1}(W - \overset{\circ}{X}, M_1; \Lambda)$ and with it $H_{q-1}(W', M_1; \Lambda)$ vanish. Furthermore $H_q(W - \overset{\circ}{X}, M_1; \Lambda)$ is isomorphic to $U_1$ and if we equip $H_q(W - \overset{\circ}{X}, M_1; \Lambda)$ with the preferred basis of $U_1$ then $d : H_{q+1}(W', W - \overset{\circ}{X}; \Lambda) \to H_q(W - \overset{\circ}{X}, M_1; \Lambda)$ is a simple isomorphism.

Since $H_*(W, M_1)$ vanishes for $* \neq q$, the boundary operator

$$H_{2q}(W, W - \overset{\circ}{X}; \Lambda) \overset{\cong}{\longrightarrow} H_{2q-1}(W - \overset{\circ}{X}, M_1; \Lambda)$$

is an isomorphism. If we equip $H_{2q-1}(W - \overset{\circ}{X}, M_1; \Lambda)$ with the preferred basis of $H_{2q}(W, W - \overset{\circ}{X}; \Lambda)$, then on the one hand the vertical exact sequence is a sequence of based modules and the Whitehead torsion of this exact sequence vanishes. On the other hand, this basis is given by $S^q \times S^{q-1}$ which is the image of the preferred basis on $H_{2q}(W', W - \overset{\circ}{X}; \Lambda)$ under $d : H_{2q}(W', W - \overset{\circ}{X}; \Lambda) \longrightarrow H_{2q-1}(W - \overset{\circ}{X}, M_1; \Lambda)$. Thus the Whitehead torsion of the vertical sequence vanishes, too.

Since $(W, M_1)$ and $(W, W - \overset{\circ}{X})$ have trivial Whitehead torsion and the Whitehead torsion of the acyclic complex given by the based horizontal homology sequence vanishes, the additivity formula implies that $(W - \overset{\circ}{X}, M_1)$ has trivial torsion. The same argument applied to the vertical sequence implies that $(W', M_1)$ has trivial torsion. Thus $W'$ is an $s$-cobordism.

Now, we show that $\theta$ is elementary if $(W, \bar{\nu})$ is bordant relative to the boundary to an $s$-cobordism $(N, \bar{\nu})$. We can suppose that $\bar{\nu}$ is a $q$-equivalence. If the bordism to an $s$-cobordism is obtained by a sequence of additions of handles on disjoint embeddings $(S^q \times D^q)_i$ as above, then the considerations above show in turn that $\theta\,(W, \overline{\nu})$ is elementary. But after possibly stabilizing $(W, \overline{\nu})$ by connected sum with $(S^q \times S^q)'s$, a bordism of this type always exists. Namely, by similar considerations (as in §4), we can transform $(N, \overline{\nu})$ by surgeries on disjoint embeddings $(S^{q-1} \times D^{q+1})_i$ into a normal $(q-1)$-smoothing $(N', \overline{\nu})$. In turn $(N, \overline{\nu})$ is obtained from $(N', \overline{\nu})$ by surgeries on $(D^q \times S^q)_i$. On the other hand, since $(N', \overline{\nu})$ and $(W, \overline{\nu})$ are bordant normal $(q-1)$-smoothings they are stably diffeomorphic by Theorem 2. This ends the proof of Theorem 3. $\qquad\square$

*Remark.* If $B$ is a finite simple Poincaré complex and the normal smoothings are simple homotopy equivalences then Theorem 3 is the same as Wall's result [W1, Th. 6.4].

*Remark.* One can also ask for obstructions for replacing $(W, \bar{\nu})$ as in Theorem 3 by an $h$-cobordism instead of an $s$-cobordism. The only difference in the proof of Theorem 3 is that one could drop the bases everywhere. Thus one would have to modify the obstruction monoids (or groups) by omitting the bases everywhere. This remark might be helpful in understanding the definition of $L_{2q}(\pi, w)$ where we require that $V$ be based but do not require that the adjoint of $\lambda$ be simple.

## 6. The main theorem for odd-dimensional $W$

We begin with the definition of the obstruction monoid $l_{2q+1}(\pi, w)$. An object is represented by a half rank based direct summand $V$ in a hyperbolic form of rank $r$ on $\Lambda^{2r}$. More precisely we consider for $\varepsilon = (-1)^q$ pairs $(H_\varepsilon^r, V)$, where $V$ is a based submodule of rank $r$ in $\Lambda^{2r}$ and $V$ is a direct summand. We stabilize these objects by identifying $(H_\varepsilon^r, V)$ with $(H_\varepsilon^r \perp H_\varepsilon, V \perp \Lambda \times \{0\})$.

Now, we use as in Wall's book [W1] the group $RU^\varepsilon(\Lambda) = \lim RU^\varepsilon(\Lambda^r)$, where $RU^\varepsilon(\Lambda^r)$ is generated by the flip $\sigma$ mapping $e \mapsto \varepsilon \mathrm{f}$ and $f \mapsto e$ and those simple isometries of $H_\varepsilon^r$ preserving $\Lambda^r \times \{0\}$ and inducing a simple isomorphism on $\Lambda^r \times \{0\}$. The group $RU^\varepsilon(\Lambda)$ acts on the stable equivalence class of objects $(H_\varepsilon^r, V)$ by mapping $V$ to $A \cdot V$ and the resulting set of equivalence classes is a monoid under orthogonal sum. Similarly, as in the definition of $l_{2q}$, we obtain - for $q$ odd - a modified monoid by passing from $\mu$ to $\tilde{\mu}$ with values in $\Lambda/\langle x - \varepsilon \bar{x}, 1 \rangle$. We denote the corresponding group of isometries $RU^\varepsilon(\Lambda)^\sim$.

*Definition.* The monoid $l_{2q+1}(\pi, w)$ is given by equivalence classes of stable pairs $(H_\varepsilon^r, V)$, where $V$ is a based submodule of rank $r$ in $\Lambda^{2r}$ and $V$ is a direct summand and where two such stable pairs are equivalent if they are in the same orbit under the action of $RU^\varepsilon(\Lambda)$. If, for $q$ odd, we take the same objects but replace $RU^\varepsilon(\Lambda)$ by $RU^\varepsilon(\Lambda)^\sim$, we obtain instead $l_{2q+1}^\sim(\pi, w)$.

The reader should note that in the definition of pairs $(H_\varepsilon^r, V)$ the preferred *lagrangian* (i.e. a based half rank direct summand on which the form vanishes identically, a subkernel in the language of [W1]) $\Lambda^r \times 0$ plays an essential role. We have not made this explicit in the notation since this information is implicitly contained in the basis of $H_\varepsilon^r$, but its role will become clear when we define elementary obstructions.

Similarly as for $l_{2q}(\pi, w)$ we consider a submonoid of $l_{2q+1}(\pi, w)$ consisting of those pairs $(H_\varepsilon^r, V)$ where $\lambda$ and $\mu$ are trivial on $V$. This is an *abelian group* denoted $L_{2q+1}(\pi, w)$. To see this, one has to describe an inverse. This is a rather delicate point. We first note that for those pairs $(H_\varepsilon^r, V)$ where $\lambda$ and $\mu$ are trivial on $V$ there is an $A \in U^\varepsilon(\Lambda)$ with $A(\Lambda^r \times \{0\}) = V$. Then we use Lemma 6.2 from [W1] which says that $A \oplus A^{-1} \in RU^\varepsilon(\Lambda)$. This lemma is only proved, in [W1], for $A \in SU^\varepsilon(\Lambda)$ but the same argument works for $A \in U^\varepsilon(\Lambda)$. Thus the inverse of $(H_\varepsilon^r, V)$ is given by $(H_\varepsilon^r, A^{-1}(\Lambda^r \times \{0\}))$.

Again we have a homomorphism $L_{2q+1}(\pi, w) \rightarrow \mathrm{Wh}(\pi)$ mapping $(H_\varepsilon^r, V)$ to the Whitehead torsion of the base change between the standard basis on $H_\varepsilon^r$ and the basis obtained from the basis $\alpha_1, \ldots, \alpha_r$ on $V$ and dual elements $\beta_1, \ldots, \beta_r$ with $\lambda(\alpha_i, \beta_j) = \delta_{i,j}$. The kernel of this homomorphism is by definition the *Wall group* $L_{2q+1}^s(\pi, w)$ [W1].

Now, let $M_0$ and $M_1$ be $(2q)$-dimensional manifolds and $f : \partial M_0 \rightarrow \partial M_1$ be a diffeomorphism. Suppose that there are normal $(q-1)$-smoothings in a fibration $B$ over $BO$ compatible with $f$. Let $W$ together with a normal $B$-structure $\bar{\nu}$ be a $B$-zero bordism of $M_0 \cup_f M_1$. Then we can by Proposition 4 assume that $\bar{\nu}$ is a $q$-equivalence; i.e., $(W, \bar{\nu})$ is a normal $(q-1)$-smoothing. From this one can easily conclude that $\pi_q(W) \rightarrow \pi_q(W, M_i)$ is surjective and that the boundary operator $\pi_{q+1}(B, W) \longrightarrow \pi_q(W, M_i)$ is also surjective. Since $(B, W)$ is $q$-connected and $B$ has finite $(q+1)$-skeleton $\pi_{q+1}(B, W) \cong H_{q+1}(B, W; \Lambda)$ is finitely generated. We choose disjoint embeddings $(S^q \times D^{q+1})_i \subset W$ compatible with the $B$-structure representing generators of $\mathrm{im}(d : \pi_{q+1}(B, W) \to \pi_q(W))$ and denote $\cup_i(S^q \times D^{q+1})_i$ by $U$.

We consider the exact sequence with $\Lambda$-coefficients:

$$H_{k+1}(W - \overset{\circ}{U}, \partial U \cup M_0) \longrightarrow H_k(\partial U) \longrightarrow H_k(W - \overset{\circ}{U}, M_0).$$

From now on we suppose that $\chi(M_0) = \chi(M_1)$, where $\chi$ denotes the Euler characteristic. Standard arguments in algebraic topology show (compare [W1, Lemma 2.3 and p. 50]) :

  i) This sequence vanishes except for $k = 0, q, 2q$. In the cases $k = 0, 2q$ the left or right maps are obviously isomorphisms and the corresponding modules are free with a canonical geometric basis. For $k = q$ all terms are stably free and after stabilization we equip these modules with a preferred basis such that the Whitehead torsion of all three pairs is 0.

 ii) $V := H_{q+1}(W - \overset{\circ}{U}, \partial U \cup M_0; \Lambda)$ is a half rank direct summand in $H_q(\partial U; \Lambda)$.

We denote $w_i(B) := p^\star w_i(BO)$ the $q^{\text{th}}$ Stiefel-Whitney class, where $p : B \to BO$ is the projection.

*Definition.* For $(W, \bar\nu)$ with the properties above and either $q \neq 3, 7$ or $q = 3, 7$ and $\langle w_{q+1}(B), \pi_{q+1}(B)\rangle = 0$, define $\theta(W, \bar\nu) \in l_{2q+1}(\pi_1(B), w_1(B))$ by:

$$\theta(W, \bar\nu) := [H_q(\partial U; \Lambda), V],$$

where $V$ is based as under i). If $q = 3, 7$ and $\langle w_{q+1}(B), \pi_{q+1}(B)\rangle \neq 0$, then we consider $\theta(W, \bar\nu) \in l^\sim_{2q+1}(\pi, w)$ instead.

We will show below that this invariant is a well-defined bordism invariant relative to the boundary. The definition of this invariant is with respect to $M_0$. Of course one could use $M_1$ instead and would get a different invariant $(H_q(\partial U; \Lambda), V')$ which carries the same information since $V' = V^\perp$ with the induced basis, as one can easily check.

*Remark.* If the smoothings on $M_i$ are normal $k$-smoothings for some $k \geq q$, the invariant is contained in the subgroup $L_{2q+1}(\pi_1(B), w_1(B))$ and if, in addition, $B$ is a finite simple Poincaré complex and the smoothings are simple homotopy equivalences the obstruction sits in $L^s_{2q+1}(\pi_1(B), w_1(B))$ and is equal to Wall's surgery obstruction [W1, 6.1].

*Definition.* We call an element $\theta \in l_{2q+1}(\pi, w)$ (resp. $\theta \in l^\sim_{2q+1}(\pi, w)$) *elementary* if and only if it has a representative $(H^r_\varepsilon, V)$ with $V \oplus \{0\} \times \Lambda^r = H^r$ and the basis of $V$ together with the standard basis of $\{0\} \times \Lambda^r$ is equivalent to the standard basis of $\Lambda^{2r}$.

The following equivalent formulation of this definition is sometimes useful. A *lagrangian complement* of a half rank direct summand $V$ in $H(\Lambda^r)$ is a based half rank direct summand $\hat{V}$ in $\Lambda^{2r}$ such that $\lambda$ and $\mu$ vanish on $\hat{V}$ and $V \oplus \hat{V} = \Lambda^{2r}$ as based modules. Now $\theta \in l_{2q+1}(\pi, w)$ (resp. $\theta \in \widetilde{l_{2q+1}}(\pi, w)$ ) is elementary, if and only if there is a representative $(H_\varepsilon^r, V)$ of $\theta$ such that $V$ and $\Lambda^r \times \{0\}$ have a common lagrangian complement.

Next, we can formulate the main result for $W$ odd-dimensional.

THEOREM 4. *Let $M_0$ and $M_1$ be connected $2q$-dimensional manifolds ($q \geq 2$) with the same Euler characteristic and $f : \partial M_0 \to \partial M_1$ be a diffeomorphism. Suppose that there are normal $(q-1)$-smoothings in a fibration $B$ over $BO$ compatible with $f$. Let $W$ together with a normal $B$-structure $\bar{\nu}$ be a $B$-zero bordism of $M_0 \cup_f M_1$. Then $\theta(W, \bar{\nu}) \in l_{2q+1}(\pi_1(B), w_1(B)$ (resp. $\theta(W, \bar{\nu}) \in \widetilde{l_{2q+1}}(\pi_1(B), w_1(B)$, if $q = 3, 7$ and $\langle w_{q+1}(B), \pi_{q+1}(B) \rangle \neq 0$) is a bordism invariant relative to the boundary and $(W, \bar{\nu})$ is bordant relative to the boundary to a relative $s$-cobordism if and only if $\theta(W, \bar{\nu})$ is elementary.*

*Proof.* Let $(W, \bar{\nu})$ be a normal $(q-1)$-smoothing. Consider a fixed set of generators $x_1, \ldots, x_k$ of $\mathrm{im}(\pi_{q+1}(B, W) \longrightarrow \pi_q(W))$ and represent them by compatible embeddings $(S^q \times D^{q+1})_i$. We first discuss the effect of a change of the framings (compatible with the $B$-structure). They correspond to composition with diffeomorphisms $S^q \times D^{q+1} \to S^q \times D^{q+1}$ mapping $(x, y)$ to $(x, \alpha(x) \cdot y)$ for an appropriate map $\alpha : S^q \to O(q+1)$. Such a composition leads again to a compatible embedding if and only if $\alpha \in \pi_q(O(q+1)) \cong \pi_{q+1}(BO(q+1))$ maps to 0 under $\pi_{q+1}(BO(q+1)) \to \pi_{q+1}(BO) \to \pi_q(F)$, where $F$ is the fibre of $B \to BO$ (see proof of Lemma 2). The induced map on $H_q(S^q \times S^q; \Lambda)$ maps $e \mapsto e + \deg(\rho\alpha) \cdot f$ and $f \mapsto f$, where $\rho$ is the evaluation map $SO(q+1) \to S^q$. For $q \neq 1, 3, 7$ the degree $\deg(\rho\alpha)$ is always even and thus the induced map in homology is contained in $RU^\varepsilon(\Lambda)$, showing that the invariant $\theta(W, \bar{\nu})$ does not depend on the choice of the compatible framing. If $q = 1, 3, 7$, there is an $\alpha$ with $\deg(\rho\alpha) = 1$. The degree $\deg(\rho\alpha)$ is mod 2 equal to $\langle w_{q+1}(E_\alpha), [S^{q+1}] \rangle$, where $E_\alpha$ is the vector bundle over $S^{q+1}$ corresponding to $\alpha$. Thus, if there is no element in $\pi_{q+1}(B)$ on which $w_{q+1}$ evaluates nontrivially, i.e. $\langle w_{q+1}(B), \pi_{q+1}(B) \rangle = 0$, then the composition with the corresponding diffeomorphism is not compatible, since $\alpha \in \pi_q(O(q+1)) \cong \pi_{q+1}(BO(q+1))$ maps nontrivially under $\pi_{q+1}(BO(q+1)) \to \pi_{q+1}(BO) \to \pi_q(F)$. Thus again only those changes of framings are possible where the induced map in homology is contained in $RU^\varepsilon(\Lambda)$ showing that the invariant $\theta(W, \bar{\nu})$ does not depend on the choice of the compatible framings if $\langle w_{q+1}(B), \pi_{q+1}(B) \rangle = 0$. On the other hand if $\langle w_{q+1}(B), \pi_{q+1}(B) \rangle \neq 0$, then there is a compatible change of framing whose induced map in homology $H_q(S^q \times S^q; \Lambda)$ maps $e \mapsto e + f$ and $f \mapsto f$. This is not contained in $RU^\varepsilon(\Lambda)$ but in $RU^\varepsilon(\Lambda)^\sim$ and so in the case $\langle w_{q+1}(B), \pi_{q+1}(B)\rangle \neq 0$, the invariant is only well-defined in $l\widetilde{2}_{q+1}(\pi_1(B), w_1(B))$. An easy consideration shows that $RU^\varepsilon(\Lambda)^\sim$ is generated by the isometries mapping $e_i \mapsto e_i + f_i$ and $f_i \mapsto f_i$ and by $RU^\varepsilon(\Lambda)$. For the rest of the argument keep in mind that we can realize the map $e_i \mapsto e_i + f_i$ and $f_i \mapsto f_i$ by the change of framing above. This is the only place where the argument for $q = 3, 7$ and $\langle w_{q+1}(B), \pi_{q+1}(B)\rangle \neq 0$ differs from the other cases.

The rest of the proof goes along the same scheme as in [W1, §6]. The embeddings $(S^q \times \{0\})_i$ are uniquely determined by $x_i$ up to regular homotopy, i.e. an immersion of $S^q \times I$ into $W$ extending the two embeddings. By the same argument as in [W1, p. 58], the invariant $\theta(W)$ changes by the action of an element in the subgroup of $RU_k^\varepsilon(\Lambda)$ given by those isometries fixing $\Lambda^k \times \{0\}$ identically. By choice of an appropriate regular homotopy each element of this subgroup occurs.

Thus the invariant is described once a set of generators of $\mathrm{im}(\pi_{q+1}(B, W) \longrightarrow \pi_q(W))$ is chosen. To show that it is independent of the set of generators, it is enough to show that it is the same for $x_1, \ldots, x_k$ and $x_1, \ldots, x_k, 0$ and for $x_1, \ldots, x_k, 0$ and $x_1, \ldots, x_k, y$ for an arbitrary $y$ in $\mathrm{im}(\pi_{q+1}(B, W) \longrightarrow \pi_q(W))$. For, then we can inductively go from one set of generators to another one, if we also allow permutation of the generators, something that does not change the equivalence class. The first step corresponds to changing the invariant by stabilizing (note that we have a slightly different convention from [W1], where Wall stabilizes instead by $\sigma$). The second step can, using permutations, be replaced by a sequence of one of the following steps: replace the first element by its product with $\pm g$ for some $g \in \pi_1(B)$ or replace the first element by the sum of the first two elements. Obviously, neither step changes our equivalence class. The first changes the obstruction by the isometry mapping $e_1 \mapsto \pm g \cdot e_1$, $f_1 \mapsto \pm w_1(B)(g) \cdot f_1$ and fixing $e_i$ and $f_i$ for $i \geq 2$. The second corresponds to changing by the isometry mapping $e_1 \mapsto e_1 + e_2$, $e_2 \mapsto e_2$, $f_1 \mapsto f_1$ and $f_2 \mapsto f_2 - f_1$ and fixing $e_i$ and $f_i$ for $i > 2$. Thus we have shown that for fixed $W$ the invariant $\theta(W, \bar{\nu})$ is independent of all choices.

Now, we show the bordism invariance. We first note that surgery on $(S^q \times D^{q+1})_i$ replaces $\theta(W, \overline{\nu})$ by the action of

$$\begin{pmatrix} 1 & & & & 0 \\ & \ddots & & & \\ & & \sigma & & \\ & & & \ddots & \\ 0 & & & & 1 \end{pmatrix}.$$

On the other hand, if $W$ and $W'$ are bordant relative to the boundary under a highly connected bordism, then one can pass from $W$ to $W'$ by a sequence of surgeries where the cores of the embeddings $S^q \times D^{q+1}$ represent classes in $\mathrm{im}(\pi_{q+1}(B, W) \longrightarrow \pi_q(W))$, and the surgeries are compatible with the $B$-structure (compare [W1, p. 61]). As we are free in the choice of our system of embeddings representing generators of $\mathrm{im}(\pi_{q+1}(B, W) \longrightarrow \pi_q(W))$ we can assume that the surgeries are all performed on $(S^q \times D^{q+1})'s$ contained in $U$, proving the bordism invariance.

Now, we want to show that if $\theta(W, \overline{\nu})$ is elementary, $W$ is bordant relative to the boundary to an $s$-cobordism. For this we first show that if we apply the action of an element of $RU^\varepsilon(\Lambda)$ to $(H_q(\partial U; \Lambda), V)$ the resulting pair is equal to $(H_q(\partial U'; \Lambda), V')$ for some $(W', \overline{\nu'})$, which has the same properties as $W$ and is bordant to $W$ relative to the boundary. The action of $RU^\varepsilon(\Lambda)$ corresponds on the one hand to stabilization which can geometrically be realized with the same $W$ by adding to $U$ an embedding $S^q \times D^{q+1}$ which is contained in a ball $D^{2q+1} \subset W$ disjoint to the other embeddings. On the other hand we have to realize the action of an element $\alpha \in RU^\varepsilon(\Lambda)$. The group $RU^\varepsilon(\Lambda)$ is generated by the following isometries: (a) the flip $\sigma$, (b) permutation of the hyperbolic summands, (c) the isometry mapping $e_1 \mapsto \pm g \cdot e_1$, $f_1 \mapsto \pm w_1(B)(g) \cdot f_1$ and fixing $e_i$ and $f_i$ for $i \geq 2$, where $g \in \pi_1(B)$, (d) by the isometry mapping $e_1 \mapsto e_1 + e_2$, $e_2 \mapsto e_2$, $f_1 \mapsto f_1$ and $f_2 \mapsto f_2 - f_1$ and finally (e) by the isometries mapping $e_i \mapsto e_i$ and $f_i \mapsto f_i + \sum c_{ij} e_j$, where $c_{ij} = \varepsilon \bar{c}_{ji}$ and $c_{ii}$ is of the form $c_i - \varepsilon \bar{c}_i$ [W1, pp. 57–60]. We have to realize all these transformations geometrically. (a) The flip interchanging $e_i$ and $f_i$ corresponds to carrying out surgery on $(S^q \times D^{q+1})_i$ so that $W$ is replaced by $W'$, which has the same properties as $W$ and is bordant to $W$ relative to the boundary. (b) Permutation corresponds to permutation of the components of $U$. (c) This transformation corresponds to the action of $\pi_1$ on generators of $\mathrm{im}(\pi_{q+1}(B, W) \longrightarrow \pi_q(W))$. (d) This corresponds to a base change of $\mathrm{im}(\pi_{q+1}(B, W) \longrightarrow \pi_q(W))$ replacing $(S^q \times D^{q+1})_1$ by the fibre bundle connected sum with $(S^q \times D^{q+1})_2$ and leaving the other embeddings unchanged. (e) These transformations form the subgroup given by those isometries fixing $\Lambda^k \times \{0\}$ identically and as mentioned above Wall shows that one can realize this action by changing the embeddings through appropriate regular homotopies.

By these arguments we can now assume that we have found a $W'$ in the $B$-bordism class relative to the boundary of $W$ such that the invariant $(H_q(\partial U'; \Lambda), V')$ has the property: $V' \oplus \{0\} \times \Lambda^r = \Lambda^{2r} = H_q(\partial U'; \Lambda)$ and the basis of $V'$ together with the standard basis of $\{0\} \times \Lambda^r$ is equivalent to the standard basis of $\Lambda^{2r}$. We claim that then $W'$ is an $s$-cobordism. For this we consider the exact sequences with $\Lambda$-coefficients (see Figure 1.)

Since all other homology groups of $(W, M_0)$ vanish and $\pi_1(M_i) \cong \pi_1(W)$, $W$ is an $s$-cobordism if and only if the map $H_{q+1}(U', \partial U'; \Lambda) \to H_q(W - \overset{\circ}{U'}, M_0; \Lambda)$ is a simple isomorphism, where $H_q(W - \overset{\circ}{U'}, M_0; \Lambda)$ is based in such a way that the vertical sequence has trivial Whitehead torsion. But

$$0$$
$$\downarrow$$
$$H_{q+1}(W - \overset{\circ}{U'}, M_0 \cup \partial U')$$
$$\downarrow$$

$$H_{q+1}(U', \partial U') \quad \longrightarrow \quad H_q(\partial U')$$
$$\downarrow \cong \qquad\qquad\qquad \downarrow$$
$$0 \to H_{q+1}(W, M_0) \to H_{q+1}(W, W - \overset{\circ}{U'}) \to H_q(W - \overset{\circ}{U'}, M_0) \quad \to \quad H_q(W, M_0) \to 0.$$
$$\downarrow$$
$$0$$

Figure 1.

$H_{q+1}(U', \partial U'; \Lambda)$ maps in $H_q(\partial U'; \Lambda)$ injectively to the based submodule with basis $f_i$, thus to $\{0\} \times \Lambda^r$ with the canonical basis. Then the fact that in the vertical sequence the image of $H_{q+1}(W - \overset{\circ}{U'}, M_0 \cup \partial U'; \Lambda)$ in $H_q(\partial U'; \Lambda)$ is $V'$ and that $V' \oplus \{0\} \times \Lambda^r = \Lambda^{2r} = H_q(\partial U'; \Lambda)$ and the basis of $V'$ together with the standard basis of $\{0\} \times \Lambda^r$ is equivalent to the standard basis of $\Lambda^{2r}$ implies the desired statement.

To finish the proof we have to show that if in turn $(W, \overline{\nu})$ is bordant relative to the boundary to an $s$-cobordism, then $\theta(W, \overline{\nu})$ is elementary. Obviously this is the case if $W$ is an $s$-cobordism and the statement follows since the invariant is a cobordism invariant. $\qquad \square$

*Remark.* If $B$ is a finite simple Poincaré complex and the normal smoothings are simple homotopy equivalences then Theorem 4 is the same as Wall's result [W1, Th. 5.6].

*Remark.* One can also ask for obstructions for replacing $(W, \bar{\nu})$ as in Theorem 4 by an $h$-cobordism instead of an $s$-cobordism. As in the even-dimensional case the only difference in the proof of Theorem 4 is that one could drop the bases everywhere. Thus one would have to modify the obstruction monoids (or groups) by omitting the bases everywhere. This remark might be helpful in understanding the definition of $L_{2q+1}(\pi, w)$ where we require that $V$ be based but do not require that the base change between the standard basis on $H_\varepsilon^r$ and the basis obtained from the basis $\alpha_1, \dots, \alpha_r$ on $V$ and dual elements $\beta_1, \dots, \beta_r$ with $\lambda(\alpha_i, \beta_j) = \delta_{i,j}$ be simple.

## 7. Analysis of obstructions in $l_{2q+1}$ under some stability assumptions

The aim of this section is to prove that in the situation of Theorem 4 the obstruction $\theta(W, \overline{\nu}) \in l_{2q+1}(\pi_1(B), w_1(B))$ is elementary if $B$ is simply connected and $q$ odd, or $q$ even and $\mathrm{Ker}(\pi_q(M_0) \to \pi_q(B))$ splits off a hyperbolic plane $H_+(\mathbb{Z})$, or if $\pi_1(B)$ is finite and $\mathrm{Ker}(\pi_q(M_0) \to \pi_q(B))$ splits off two hyperbolic planes: $H_{(-1)^q}(\Lambda^2)$. We only discuss obstructions in $l_{2q+1}$ and leave the obvious modifications for the case of obstructions in $\widetilde{l_{2q+1}}$ to the reader. Denote $\varepsilon = (-1)^q$. Recall that a lagrangian complement of a half rank direct summand $V$ in $H_\varepsilon(\Lambda^r)$ is a based half rank direct summand $\hat{V}$ in $\Lambda^{2r}$ such that $\lambda$ and $\mu$ vanish on $\hat{V}$ and $V \oplus \hat{V} = \Lambda^{2r}$ as based modules. We begin with the following proposition:

PROPOSITION 8. *For $W$ as in Theorem 4 write $\theta(W, \overline{\nu}) = (H_\varepsilon(\Lambda^r), V)$.*

i) *There is a surjective isometry of quadratic forms*

$$V \longrightarrow Ker(\pi_q(M_0) \longrightarrow \pi_q(B)).$$

ii) *Let $\pi : S \longrightarrow Ker(\pi_q(M_0 \longrightarrow \pi_q(B))/\mathrm{rad}$ be a surjection from a free based $\Lambda$-module, where $\mathrm{rad}$ is the radical consisting of all $x$ with $\lambda(x, y) = 0$ for all $y$ and $\mu(x) = 0$. Equip $S$ via $\pi$ with a quadratic form. Then there is an isometric embedding of $S$ into a half rank direct summand of $H_\varepsilon(\Lambda^s)$ and $\rho \in L_{2q+1}(\pi_1(B), w_1(B))$ such that*

$$\theta(W, \overline{\nu}) = [H_\varepsilon(\Lambda^s), S] + \rho.$$

iii) *If $(H_\varepsilon(\Lambda^r), V)$ has a lagrangian complement, then there is a $\rho \in L_{2q+1}(\pi_1, w_1)$ such that $(H_\varepsilon(\Lambda^r), V) \perp \rho$ is elementary.*

*Proof.* i) Consider the sum of boundary operators

$$H_{q+1}(W - \overset{\circ}{U}, \partial U \cup M_0; \Lambda) \to H_q(\partial U, \Lambda) \oplus \mathrm{Ker}(\pi_q(M_0) \to \pi_q(B)).$$

The first component $d_1$ of this map maps $H_{q+1}(W - \overset{\circ}{U}, \partial U \cup M_0; \Lambda)$ isomorphically to $V$ and the second component $d_2$ is surjective since $U$ generates $\mathrm{Ker}(\pi_q(W) \to \pi_q(B))$. By Proposition 5, $\lambda$ and $\mu$ vanish on image $(d_1, d_2)$. Thus $d_2 d_1^{-1} : V \longrightarrow \mathrm{Ker}(\pi_q(M_0) \to \pi_q(B))$ is a surjective isometry by the fact that, since $\partial U$ and $M_0$ are disjoint, intersection numbers between cycles in them are zero.

ii) Suppose $\theta(W, \overline{\nu})$ is $(H_\varepsilon(\Lambda^r), V)$ with isometric projection

$$p : V \longrightarrow K := \mathrm{Ker}\,(\pi_q(M_0) \to \pi_q(B))/\mathrm{rad}$$

as in i). After perhaps stabilizing $\theta(W, \overline{\nu})$ by orthogonal sum with $(H_\varepsilon(\Lambda^s), \Lambda^s \times \{0\})$ we construct a surjective homomorphism $V \longrightarrow S$ commuting with the projections $p$ and $\pi$. From this, one obtains, after adding to $V$ the module $S$ with the trivial map to $S$ and using a splitting of $V \oplus S \to S$, a commutative diagram

$$
\begin{array}{ccc}
V \oplus S & \xrightarrow[\Phi]{\cong} & S \oplus V \\
\searrow p+0 & & \swarrow \pi+0 \\
& K &
\end{array} \quad .
$$

Pull the quadratic form on $K$ back via $p+0$ and $\pi+0$ to obtain quadratic forms such that $\Phi$ is an isometry. After perhaps stabilizing $V$ further and composing with an appropriate automorphism of $V$ on the right side, one can assume that $\Phi$ is a simple isometry.

Choose a simple isomorphism $S \to \Lambda^s$ and identify the trivial element in $l_{2q+1}(\pi_1(B), w_1(B))$ given by $(H_\varepsilon(\Lambda^s), \Lambda^s \times \{0\})$ with $(H_\varepsilon(S), S \times \{0\})$. Then

$$
\begin{aligned}
\theta(W, \overline{\nu}) &= [H_\varepsilon(\Lambda^r), V] = [H_\varepsilon(\Lambda^r), V] \perp [H_\varepsilon, (V_1), S \times \{0\}] \\
&= [H_\varepsilon(\Lambda^{r+s}, V \oplus S] = [H_\varepsilon(\Lambda^{r+s}), \Phi^{-1}(S \oplus V)],
\end{aligned}
$$

where on the right side the form on $S \oplus V$ is induced by $\pi+0$. In particular, $S$ is orthogonal to $V$.

Since the quadratic form vanishes on $\Phi^{-1}(V)$ and this is a direct summand in $\Lambda^{2(r+s)}$, we can embed $\Phi^{-1}(V)$ into a lagrangian and, since the group of isometries of a hyperbolic form acts transitively on the lagrangians, there is an isometry $A$ such that $A\Phi^{-1}(V) = \langle e_1, \ldots, e_r \rangle$. After perhaps embedding $S$ differently into $S \oplus V$ by an isometry we can assume $A\Phi^{-1}(S) \subset \langle e_1, \ldots, e_r, ; f_1, \ldots, f_r \rangle^\perp$. Thus,

$$
[H_\varepsilon(\Lambda^{r+s}), A(\Phi^{-1}(S \oplus V))] = [H_\varepsilon(\Lambda^s), A\Phi^{-1}(S)].
$$

On the other hand, if $B$ is an isometry of $H(\Lambda^t)$ and $U \subset H_\varepsilon(\Lambda^t)$ is a based half rank direct summand, then

$$
[H_\varepsilon(\Lambda^t), B(U)] = [H_\varepsilon(\Lambda^{2t}), (\mathrm{Id} \oplus B)(\Lambda^t \times \{0\} \oplus \{0\} \times U)]
$$
$$
[H_\varepsilon(\Lambda^{2t}), (B \oplus \mathrm{Id})(\Lambda^t \times \{0\} \oplus \{0\} \times U)]
$$
$$
= [H_\varepsilon(\Lambda^t), B(\Lambda^t \times \{0\})] \perp [H_\varepsilon(\Lambda^r), U]
$$

and with this we obtain

$$
[H_\varepsilon(\Lambda^{r+s}), A(\Phi^{-1}(S \oplus V))]
$$
$$
= [H_\varepsilon(\Lambda^{r+s}), A(\Lambda^{r+s} \times \{0\})] \perp [H_\varepsilon(\Lambda^{r+s}), \Phi^{-1}(S \oplus V)]
$$

and the first summand is in $L_{2q+1}(\pi_1(B), w_1(B))$.

iii) Let $\hat{V}$ be a lagrangian complement of $V$. Then there is an isometry $A$ mapping $V$ to $\{0\} \times \Lambda^r$. Thus $[H_\varepsilon(\Lambda^r), A(V)]$ is elementary and $[H_\varepsilon(\Lambda^r), A(V)] = [H_\varepsilon(\Lambda^r), V] \perp [H_\varepsilon(\Lambda^r, A(\Lambda^r \times \{0\}))] = [H_\varepsilon(\Lambda^r), V] \perp \rho$ with $\rho \in L_{2q+1}(\pi_1, w_1)$. $\qquad \square$

Now, we prove the theorem announced at the beginning of this section.

THEOREM 5. *For $W$ as in Theorem 4 there is another bordism $W'$ between the same ends such that $\theta(W', \overline{\nu})$ is elementary if one of the following conditions is fulfilled.*

i) *$q$ odd and $B$ 1-connected,*

ii) *$q$ even, $B$ 1-connected and $\mathrm{Ker}\,(\pi_q(M_0) \to \pi_q(B))/_{\mathrm{rad}}$ splits off $H_+(\mathbb{Z})$.*

iii) *$\pi_1(B)$ is finite and $\mathrm{Ker}\,(\pi_q(M_0) \to \pi_q(B))/_{\mathrm{rad}}$ splits off $H_\varepsilon(\Lambda^2)$.*

COROLLARY 4. *Let $M_0$ and $M_1$ be $2q$-dimensional manifolds, with the same Euler characteristic, which admit $B$-bordant normal $(q-1)$-smoothings in a fibration $B$ over $BO$. Then they are diffeomorphic (homeomorphic if $q = 2$) if one of the conditions* i), ii) *or* iii) *is fulfilled.*

The statement under condition iii) was, for $q = 2$, proved in [H-K1]. A similar argument holds in higher dimensions.

*Proof.* The main ingredient is the following proposition which is a consequence of a result by Bass [Ba2].

PROPOSITION 9. *Let $V$ be a submodule of $H_\varepsilon(\Lambda^r)$ and let $H_1$ be equal to $H_\varepsilon(\Lambda)$. If either $q$ is odd and $\pi_1 = \{0\}$, or $q$ is even, $\pi_1 = \{0\}$ and $V = V' \perp H_+(\mathbb{Z})$ or $\pi_1$ is finite and $V \cong V' \perp H_\varepsilon(\Lambda^2)$ then for each hyperbolic plane $H \subset V \perp H_1$ there is an element $A$ of $RU^\varepsilon(\Lambda)$ and so an isometry (after stabilization if required) of $H_\varepsilon(\Lambda^r) \perp H_1$ mapping $H$ to $H_1$ and $V \perp H_1$ to itself.*

*Proof.* In the case of finite nontrivial $\pi_1$ the result follows immediately from [Ba2, Cor. 3.5, p. 236] when we note that the group $G_1$ (in this corollary), which acts transitively on the hyperbolic planes, is contained in $RU^\varepsilon(\Lambda)$ and preserves $V \perp H_1$. In the case of $\pi_1$ trivial one uses the fact that $RU^\varepsilon(\Lambda)$ is the full group of isometries [W1, Th. 13A.1]. Thus one only has to find some isometry mapping $H$ to $H_1$. This can be done using transvections [Ba2, p. 91] which are isometries given by two elements $u$ and $v$ with $\langle u, v \rangle = 0$, $\langle u, u \rangle = 0$ and $\langle v, v \rangle = 0$:

$$\sigma_{u,v}(x) := x + \langle v, x \rangle u - (-1)^q \langle u, x \rangle v.$$

We begin with the case $q$ odd. If $H$ and $H_1$ are equipped with standard bases $e$, $f$ and $e_1$, $f_1$ write $e = x + ae_1 + bf_1$ with $x \in V$. Since the group of isometries of $H_1$ is $\mathrm{SL}_2(\mathbb{Z})$, we can assume $b = 0$: $e = x + ae_1$. Write $f = y + be_1 + cf_1$. Now, $\sigma^r_{e,e_1}(f) = y' + (b + rca + r)e_1 + cf_1$. Thus we can assume $(b, c) = 1$ and after applying an appropriate isometry of $H_1$ obtain $f = y + e_1$. Applying $\sigma^{-1}_{y,-f_1}$ maps $f$ to $e_1$. Thus we can assume $f = e_1$ and write $e = x + ae_1 + bf_1$. Since $\langle e, f \rangle = 1$, we have $e = x + ae_1 - f_1$. Replacing $e$ by $e - af$ and making another base change we can assume $e = e_1$, $f = y + f_1$. Applying $\sigma^{-1}_{y,e_1}$ maps $e$ to $e_1$ and $f$ to $f_1$, finishing the case $q$ odd.

In case $q$ is even write $V = V' \perp H_2$ and equip $H_2$ with the standard basis $e_2, f_2$. Since the group of isometries of $H_1 \perp H_2$ acts transitively on unimodular elements of fixed length [Ba2, Th. 3.4] we can assume $e = x + k(e_2 + af_2)$ with $x \in V'$.

Since $e$ is a unimodular element in $V' \perp H_1$ there is $z \in V' \perp H_2$ with $\langle e, z \rangle = 1$ and $\langle z, z \rangle = 0$. Then $\sigma_{z, -e_1}(e) = e + e_1$. Using again the transitive actions of isometries of $H_1 \perp H_2$ on unimodular elements of fixed length, we can assume $e = x + e_1 + e_2 + rf_2$. Application of $\sigma^{-1}_{x + e_2 + rf_2, f_1}$ maps $e$ to $e_1$. Thus we can assume $e = e_1$ and write $f = y + ae_1 + f_1 + be_2 + cf_2$. Now consider the isometry which is the identity on $V'$ and on $H_1 \perp H_2$ maps $ae_1 + f_1 + be_2 + cf_2$ to $f_1 + (bc + a)e_2 + f_2$, $e_1$ to $e_1$, $e_2 \rightarrow e_2 + (c - 1)e_1$ and $f_2$ to $f_2 + (b - a - bc)e_1$. After applying this isometry we can assume $e = e_1$ and $f = y + f_1 + (bc + a)e_2 + f_2$. We finish by applying $\sigma^{-1}_{y + (bc + a)e_2 + f_2, f_1}$, mapping $e$ to $e_1$ and $f$ to $f_1 + ae_1$.  □

Now, we use this proposition to finish the proof of Theorem 5. By Theorem 2 there are a $k$ and a relative $s$-cobordism between $M_0 \sharp k (S^q \times S^q)$ and $M_1 \sharp k(S^q \times S^q)$.

After perhaps modifying $W$ by disjoint union with a closed $(2q + 1)$-dimensional $B$-manifold $X$ to obtain $W'$ we can assume that the relative $B$-bordism $W' \cup k (S^q \times D^{q+1} \times I)$ glued along $k$ disjoint embeddings of $I \times D^{2q}$ with $\{0\} \times D^{2q} \subset M_0$ and $\{1\} \times D^{2q} \subset M_1$ is $B$-bordant relative to the boundary to the $s$-cobordism. Thus $\theta$ $(W' \cup k (S^q \times D^{q+1} \times I), \bar{\nu})$ is elementary. But $\theta(W' \cup k(S^q \times D^{q+1} \times I), \bar{\nu}') = \theta(W', \bar{\nu}') \perp (H_\varepsilon(\Lambda^{2k}), H_\varepsilon(\Lambda^k))$, where the embedding of $H_\varepsilon(\Lambda^k)$ into $H_\varepsilon(\Lambda^{2k})$ is the orthogonal sum of isometric embeddings $H_\varepsilon(\Lambda)$ into $H_\varepsilon(\Lambda) \perp H_\varepsilon(\Lambda)$. We denote $H_\varepsilon(\Lambda)$ in $H_\varepsilon(\Lambda) \perp H_\varepsilon(\Lambda)$ by $H_1$ and split $H_\varepsilon(\Lambda) \perp H_\varepsilon(\Lambda) = H_1 \perp H_2$.

With the proposition above we first show inductively that $\theta(W', \bar{\nu}') = (H_\varepsilon(\Lambda^k), V)$ has a lagrangian complement. From Proposition 8, ii) we know that $V$ fulfills the assumptions of Proposition 9 above. Assume that $\theta(W', \bar{\nu}') \perp (H_1 \perp H_2, H_1)$ has a lagrangian complement. Let $\hat{V}$ be a hamiltonian complement of $V \perp H_1$. Let $e_i, f_i$ be a symplectic basis of $H_i$. Write $e_2 = x_e + y_e$ with $x_e \in V \perp H_1$ and $y_e \in \hat{V}$ and $f_2 = x_f + y_f$ with $x_f \in V \perp H_1$ and $y_f \in \hat{V}$. Then $x_e$ and $-x_f$ are the symplectic basis of a hyperbolic plane $H$ in $V \perp H_1$. Namely, using $\lambda(y_e, y_e) = 0$, since $\hat{V}$ is a lagrangian, we conclude:

$$0 = \lambda(x_e, y_e) = \lambda(x_e, y_e) + \lambda(x_e, x_e) = \lambda(x_e, y_e) + \lambda(e_2 - y_e, e_2 - y_e)$$
$$= \lambda(x_e, y_e) - \lambda(e_2, y_e) - \lambda(y_e, e_2) = \lambda(x_e, y_e) - \lambda(x_e, y_e) - \lambda(y_e, x_e)$$
$$= -\lambda(y_e, x_e).$$

Again using $0 = \lambda(x_e, y_e) + \lambda(x_e, x_e)$, we conclude: $\lambda(x_e, x_e) = 0$. Similarly $\lambda(x_f, x_f) = 0$. Using $\lambda(y_e, y_f) = 0$ and $0 = \lambda(x_f, e_2) = \lambda(x_f, y_e) + \lambda(x_f, x_e)$ and $0 = \lambda(x_e, f_2) = \lambda(x_e, y_f) + \lambda(x_e, x_f)$ we conclude:

$$1 = \lambda(e_2, f_2) = \lambda(y_e, x_f) + \lambda(x_e, y_f) + \lambda(x_e, x_f)$$
$$= -\lambda(x_e, x_f) - \lambda(x_e, x_f) + \lambda(x_e, x_f) = -\lambda(x_e, x_f).$$

Now, let $A$ be an isometry as in Proposition 9 above. Then $\hat{V}'$ $:= (A \perp Id)(\hat{V})$ is another lagrangian complement of $V \perp H_1$. Denote $B$ $:= \ker(\hat{V}' \to H_1)$ and check that $\hat{V}' = B \perp \langle e_1 - e_2, f_1 - f_2 \rangle$. Thus $B$ is a lagrangian complement of $V$.

By Proposition 8, iii) and the vanishing of $L_{2q+1}(\{e\})$ [W1] the proof that $\pi_1$ is trivial is finished. In case $\pi_1$ is nontrivial finite, the existence of a lagrangian complement implies that there is an isometry $A$ such that $[H_\varepsilon(\Lambda^r), A(V)]$ is elementary. Under the assumption that $V = V' \perp H$, where $H \cong H_\varepsilon(\Lambda^2)$, we will show that there is an isometry $B \in RU^\varepsilon(\Lambda)$ such that $BA_{|H^\perp} = \mathrm{Id}$. In particular $BA(V) = V$ and so $[H_\varepsilon(\Lambda^r), V]$ is elementary.

The existence of $B$ follows again from [Ba2, Th. 3.5, p. 236]. That is, we choose a symplectic basis $e_1, f_1, \dots, e_{r-2}, f_{r-2}$ of $H^\perp$ and apply this result inductively to find a $B \in RU^\varepsilon(\Lambda)$ such that $BA(e_i) = e_i$ and $BA(f_i) = f_i$. $\quad\square$

Now, we prove Theorem E. Recall the surgery sequence:

$$[\Sigma(M), G/O] \to L^s_{2q+1}(\pi_1(M), w_1(M)) \to S(M) \to [M, G/O]$$

where $S(M)$ is the set of ismorphism classes of pairs $f : N \to M$, $f$ a (local) orientation-preserving simple homotopy equivalence [W1]. The quotient $\pi_0(\mathrm{Aut}^s(M))/\pi_0(\mathrm{Diff}(M))$ embeds into $S(M)$ under the obvious map. Thus we consider $\pi_0(\mathrm{Aut}^s(M))/\pi_0(\mathrm{Diff}(M))$ as a subset of $S(M)$. This subset is preserved by the action of the $L$-group. Namely, the action on $S(M)$ assigns to $\theta \in L^s_{2q+1}(\pi_1(M), w_1(M))$ and a simple homotopy equivalence from $M \to M$ a simple homotopy equivalence $f : N \to M$ which is over $M$ normally bordant to the given simple homotopy equivalence $M \to M$. Since $M$ fulfills the properties of Theorem 5, there is a diffeomorphism $g : M \to N$. Thus $f : N \to M$ is in $S(M)$ equal to $f \cdot g : M \to M$. Thus we obtain the exact sequence of Theorem E.

We finish this section with a proof of Theorem F. The normal 1-type of a topological spin 4-manifold $M$ with fundamental group $\mathbb{Z}$ is $B = S^1 \times B\,\mathrm{TOPSpin} \longrightarrow B\,\mathrm{TOP}$ (we work here in the topological category and thus we replace $BO$ by $B\,\mathrm{TOP}$ and BSpin by the 3-connected cover $B\,\mathrm{TOPSpin}$ over $B\,\mathrm{TOP}$). A normal 1-smoothing is determined by the choice of a spin-structure and a generator of $\pi_1(M)$. The bordism group $\Omega_4(B)$ is equal to $\Omega_4^{\mathrm{TOPSpin}}(S^1) \cong \Omega_4^{\mathrm{TOPSpin}} \oplus \Omega_3^{\mathrm{TOPSpin}} \cong \mathbb{Z}$ determined by the signature [F-Q]. Thus any two closed topological Spin 4-manifolds with the same signature are $B$-bordant and the obstruction for finding a topological $s$-cobordism sits in $l_5(\mathbb{Z})$ (Theorem 4). Since the signature is determined by the intersection form on $\pi_2(M)$, two such manifolds with isometric intersection forms are B-bordant. To analyse the obstructions in $l_5(\mathbb{Z})$ we first note that $\pi_2(M)$ is a free $\Lambda = \mathbb{Z}[\pi_1(M)]$-module. For this ring $\Lambda$, stably free modules are actually free [Ba1]. Thus, it is enough to show that $\pi_2(M)$ is stably free. Now we note that replacing $M$ by $M\#\mathbb{C}P^2$ changes $\pi_2(M)$ by adding $\Lambda$. Thus, it is enough to show that $\pi_2(M\#\mathbb{C}P^2)$ is stably free. The normal 1-type of $M\#\mathbb{C}P^2$ is $B' = S^1 \times B\,\mathrm{STOP} \to B\,\mathrm{TOP}$, where $B\,\mathrm{STOP}$ is classifying space for stable topological oriented vector bundles. $\Omega_4(B') = \Omega_4^{\mathrm{TOP}}(S^1) \cong \mathbb{Z} \times \mathbb{Z}/2$ where the isomorphism is given by the signature and the Kirby-Siebenmann obstruction [F-Q]. Thus, $M\#\mathbb{C}P^2$ is stably homeomorphic to $S^1 \times S^3 \# r\mathbb{C}P^2$, where $r$ is the signature of $M\#\mathbb{C}P^2$ (Theorem C). Since $\pi_2(S^1 \times S^3 \# r\mathbb{C}P^2)$ is free, $\pi_2(M)$ is stably free and thus free.

Now we apply Proposition 8, ii). Modulo the sum with an element of $L_5(\mathbb{Z})$, the obstruction $\theta(W, \overline{\nu})$ for a B-bordism $W$ between two such manifolds $M_0$ and $M_1$ with the same intersection form on $\pi_2$ is given by an isometric embedding of $(\pi_2(M_0), \lambda, \mu)$ into $H_\varepsilon(\Lambda^r)$, where $r = \mathrm{rank}\ \pi_2(M_0)$. If we replace $M_0$ by $M_1$ in the definition of $\theta(W, \overline{\nu})$ we obtain an embedding of $(\pi_2(M_1), -\lambda, -\mu)$ into $H_\varepsilon(\Lambda^r)$ which is the orthogonal complement of $(\pi_2(M_0), \lambda, \mu)$ in $H_\varepsilon(\Lambda^r)$. Thus, we have an isometric embedding of $(\pi_2(M_0), \lambda, \mu) \perp (\pi_2(M_1), -\lambda, -\mu)$ into $H_\varepsilon(\Lambda^r)$. Now $\pi_2(M_0, \lambda, \mu)$ and $(\pi_2(M_1), \lambda, \mu)$ are isometric and this form is unimodular, which implies that

$$H_\varepsilon(\Lambda^r) \cong (\pi_2(M_0), \lambda, \mu) \perp (\pi_2(M_0), -\lambda, -\mu).$$

Then the diagonal embedding of $\pi_2(M_0)$ is a lagrangian complement, and we obtain from Proposition 8, iii) that $\theta(W, \overline{\nu})$ is elementary modulo the sum with an element of $L_5(\mathbb{Z})$. The Wall obstruction group $L_5(\mathbb{Z})$ is isomorphic to $\mathbb{Z}$ [W1, Th. 13A.8]. Replacing $W$ by the connected sum of $W$ with $S^1 \times r \cdot M(E_8)$, where $M(E_8)$ is the closed Spin 4-manifold with signature 8 [F-Q], we can for appropriate $r$ modify $\theta(W, \bar{\nu})$ by an arbitrary element of $L_5(\mathbb{Z})$. Thus we can assume that $\theta(W, \bar{\nu})$ is elementary and by Theorem 4 this finishes the proof of Theorem F.

*Remark.* Working with the normal 2-type instead of the normal 1-type, one can show for all oriented closed 4-manifolds with infinite cyclic fundamental group that any isometry between the intersection forms on $\pi_2$ can be realized by a homeomorphism provided that in the nonspin case the Kirby-Siebenmann obstructions agree.

Furthermore, one can classify the pseudoisotopy classes of these homeomorphisms. This problem was studied in [F-Q] and[S-W] but the answers there are slightly incorrect. An analysis using the methods developed in this paper will be given in [Kr-Te].

## 8. Applications to complete intersections

In this section we want to prove Theorem A. We begin with a reformulation. If $i : X_\delta^n \to \mathbb{C}P^\infty$ is the inclusion, the normal bundle of $X_\delta^n$ is $i^\star \xi(n, \delta)$ where $\xi(n, \delta) = -(n + r + 1) \cdot H \oplus H^{d_1} \oplus \cdots \oplus H^{d_r}$. Since $i^\star$ is injective in cohomology up to dimension $2n$ and is determined by the total degree $d$, the Pontrjagin classes of two complete intersections $X_\delta^n$ and $X_{\delta'}^n$ with the same total degree are equal if and only if $p_i(\xi(n, \delta)) = p_i(\xi(n, \delta'))$ for $2i \le n$.

Our first step is to note that if $p_i(\xi(n, \delta)) = p_i(\xi(n, \delta'))$, the two complete intersections have the same normal $(n - 1)$-type. In Proposition 3 we determined the normal $(n - 1)$-type of a complete intersection $X_\delta^n$ as $B = \mathbb{C}P^\infty \times BO\langle n + 1 \rangle @ > \xi(n, \delta) \oplus p >> BO$. Two such normal $(n - 1)$-types are equal, if the restrictions of $\xi(n, \delta)$ to $\mathbb{C}P^{[n/2]}$, the n-skeleton of $\mathbb{C}P^\infty \times BO\langle n + 1 \rangle$ are stably isomorphic. For $m \ne 2, 3 \bmod 8$ two stable real bundles over $\mathbb{C}P^m$ is stably isomorphic if and only if they have same Pontrjagin classes $p_1, \dots, p_{[m/4]}$ [Sa]. Since we also control $p_{[n/4]+1}$, the equality of the Pontrjagin classes $p_i(\xi(n, \delta))$ for $2i \le n$ implies that the restriction of $\xi(n, \delta)$ to $\mathbb{C}\mathbb{P}^{[n/2]}$, the n-skeleton of $\mathbb{C}P^\infty \times BO\langle n + 1 \rangle$ is stably isomorphic.

If $n$ is odd we have by Corollary 4 only to show that $X_\delta^n$ and $X_{\delta'}^n$ admit $B$-bordant normal $(n - 1)$-smoothings. If $n$ is even we can apply this corollary only if $\mathrm{Ker}\,(\pi_n(X_\delta^n) \to \pi_n(B))/_{\mathrm{rad}}$ splits off a hyperbolic plane. Denote the Poincaré dual of $x^{n/2}$ by $h \in H_n(X_\delta^n; \mathbb{Z})$. Then by the Hurewicz theorem $\mathrm{Ker}\,(\pi_n(X_\delta^n) \to \pi_n(B))/_{\mathrm{rad}} = h^\perp$. If $b_n(X_\delta^n) - |\mathrm{sign}(X_\delta^n)| \ge 4$ one can find a hyperbolic plane in $h^\perp$. But for $n \ge 3$ the only complete intersections with $b_{n/2}(X_\delta^n) - |\mathrm{sign}(X_\delta^n)| \le 4$ are $X_1^n$, $X_2^n$, and $X_{(2,2)}^n$ [L-W2]. We summarize these considerations as:

PROPOSITION 10. *Two complete intersections $X_\delta^n$ and $X_{\delta'}^n$ of complex dimension $n > 2$ are diffeomorphic if and only if the total degrees, the Pontrjagin classes and the Euler characteristics agree and they admit bordant normal $(n - 1)$-smoothings in $B$.*

The next step is to show that the total degree and the Pontrjagin classes determine the element in $\Omega_{2n}(B) \otimes \mathbb{Q}$. This is a standard application of the collapsing of the Atiyah-Hirzebruch spectral sequence over $\mathbb{Q}$.

Thus the difference of two complete intersections $X_\delta^n$ and $X_{\delta'}^n$ of complex dimension $n > 2$ with the same total degrees and Pontrjagin classes equipped with appropriate normal $(n - 1)$-smoothings is a torsion element in $\Omega_{2n}(B)$. Using the Pontrjagin-Thom construction we identify this group with $\pi_{2n}(M\xi(n, \delta) \wedge MO\langle n + 1 \rangle)$, where $M\xi(n, \delta)$ is the Thom spectrum of the bundle $\xi(n, \delta)$ and $MO\langle n + 1 \rangle$ is the Thom spectrum of the pullback of the universal bundle over the $n$-connected cover of $BO$. Now, one can use the Adams spectral sequence to analyze these groups. The key point is that there is a vanishing line for the Adams spectral sequence, meaning that if a torsion element has sufficiently high Adams filtration, then it is actually trivial. The mod $p$ Adams filtration of the image of an element under a map inducing the trivial map in $\mathbb{Z}/p$ homology increases by at least one. Using the inductive construction of complete intersections one can determine an upper bound of the mod $p$ Adams filtration depending on how many powers of $p$ divide the total degree. Combining this with the vanishing line leads to Theorem A. The details of this idea were carried out by Claudia Traving following suggestions by Stephan Stolz. In the following we will discuss this in more detail.

Write the total degree $d = \prod_{p \text{ prime}} p^{\nu_p(d)}$.

PROPOSITION 11 [Tr]. *The mod $p$ Adams filtration of $X_\delta^n$ equipped with an appropriate normal $(n-1)$-smoothing is at least $\nu_p(d)$.*

One can actually show that this is the precise filtration which in any case is independent of the normal $(n-1)$-smoothing since these smoothings only differ by the action of $\mathrm{Aut}(B)$ which preserves the filtration.

The proof uses an obvious translation of the process of taking the transverse intersection of manifolds into stable homotopy via Pontrjagin-Thom construction. We first note that a complete intersection actually admits a normal $(n-1)$-smoothing in the fibration $\xi(n, \delta) : \mathbb{C}P^\infty = \mathbb{C}P^\infty \times BO\langle\infty\rangle \to BO$, since the normal bundle is a pullback from a bundle over $\mathbb{C}P^\infty$. From this structure we obtain a normal $(n-1)$-smoothing over $B$ by factorization of $\xi(n, \delta) : \mathbb{C}P^\infty \to BO$ over $B$. Using this structure we consider $X_\delta^n$ via Pontrjagin-Thom construction as an element of $M\xi(n, \delta)$. Consider the map $f_{d_r} : M\xi(n, (d_1, \ldots, d_{r-1})) \to \Sigma^2(M\xi(n+1, (d_1, \ldots, d_r)))$ induced by inclusion

$$-(n+r+1)H \oplus H^{d_1} \oplus \cdots \oplus H^{d_{r-1}} \to -(n+r+1)H \oplus H^{d_1} \oplus \cdots \oplus H^{d_r}$$

(note that we assume that the Thom class of a Thom spectrum sits in dimension 0 explaining the occurrence of $\Sigma^2$). This maps the element corresponding to $X_{(d_1, \ldots, d_{r-1})}^{n+1}$ to the element corresponding to $X_{(d_1, \ldots, d_r)}^n$. Thus we are finished if $f_{d_r}$ increases the mod $p$ Adams filtration at least by $\nu_p(d_r)$. For this, one factors $f_{d_r}$ further. Write $d_r = \prod_{1 \le i \le t} s_i$ and consider the map

$$g_j : M\xi(n, d_1, \ldots, d_{r-1}, \prod_{i \le j-1} s_i) \to M\xi(n, d_1, \ldots, d_{r-1}, \prod_{i \le j} s_i)$$

induced by the $s_j$-fold tensor product mapping $H^{\prod_{i \le j-1} s_i}$ to $H^{\prod_{i \le j} s_i}$. If $p | s_j$ this map is trivial in mod $p$-homology and thus increases the Adams filtration. On the other hand we can write $f_{d_r} = g_t \cdot \ldots \cdot g_1 \cdot f_1$ and thus $f_{d_r}$ increases the mod $p$ Adams filtration at least by $\nu_p(d_r)$.

The final step is to show that torsion elements of sufficiently high filtration vanish.

PROPOSITION 12 [Tr]. *Let $X_\delta^n$ and $X_{\delta'}^n$ be complete intersections with $n \geq 3$, the same total degree $d$ and equal Pontrjagin classes. If $\nu_p(d) \geq \frac{2n+1}{2(p-1)} + 1$ for all $p$ with $p(p-1) \leq n+1$, then $X_\delta^n$ and $X_{\delta'}^n$ are, with respect to appropriate normal $(n-1)$ smoothings, $B$-bordant.*

The proof of this result in stable homotopy theory is a bit technical. Thus we will only give a sketch from which an expert should be able to fill in the details. For a spectrum $S$ with only finitely many nontrivial homotopy groups in negative dimensions and finitely generated integral cohomology in all dimensions, the mod $p$ Adams spectral sequence has the structure of a $\mathbb{Z}/p[h_0]$-module where, if $g$ is a degree $p$ map in the sphere spectrum, $h_0$ is the corresponding element in $\mathrm{Ext}_A^{1,1}(\mathbb{Z}/p, \mathbb{Z}/p) \cong \mathbb{Z}/p$, $A$ the mod $p$ Steenrod algebra. Then denote $TE_r := \{x \in E_r | h_0^n(x) = 0 \text{ for some } n \in \mathbb{N}\}$. Then, as $\mathbb{Z}/p[h_0]$-module $E_r$ splits into $TE_r$ and a free $\mathbb{Z}/p[h_0]$-module $FE_r$ and the filtration quotients of the $p$-torsion in $\pi_*(S)$ correspond to $TE_\infty$. Thus one wants to know when the map $ZTE_2^{s,t} := ZE_2^{s,t} \cap TE_2^{s,t} \to TE_\infty^{s,t}$ is surjective, where $ZE_2^{s,t}$ is the subgroup of permanent cycles in the $E_2$-term. If we from now on consider our relevant spectrum $S = M\xi(n,\delta) \wedge MO\langle n+1 \rangle$, then using information from [M-M] and some information from [Gi] about the cohomology of $BO\langle n+1 \rangle$ and the Bockstein spectral sequence for the extension $0 \to \mathbb{Z} \to \mathbb{Z} \to \mathbb{Z}/p \to 0$, one shows that for $t - s \leq 2n+1$ the map $ZTE_2^{s,t} \to TE_\infty^{s,t}$ is surjective.

Combining information from [M-M] with a vanishing result for $A_0$-free $A$-modules ($A_0$ the sub-Hopf algebra of $A$ generated by the Bockstein homomorphism $\beta$ for the extension $0 \to \mathbb{Z} \to \mathbb{Z} \to \mathbb{Z}/p \to 0$) by [A] for $p = 2$ and by [Li] for $p > 2$, one shows that for $(-1)$-connected spectra $S$ the groups $T\mathrm{Ext}_2^{s,t}$ vanish, if $s \geq 2$ and $t - s \leq 2(p-1) \cdot s - 1$. Our spectrum $S = M\xi(n,\delta) \wedge MO\langle n+1 \rangle$ is $(-1)$-connected and thus this vanishing result can be applied. The difference of our two complete intersections $X_\delta^n$ and $X_{\delta'}^n$ in $\pi_{2n}(S)$ has by assumption Adams filtration $s \geq \frac{2n+1}{2(p-1)} + 1$ for all $p$ with $p(p-1) \leq n+1$. The vanishing result implies that for $p$ with $p(p-1) \leq n+1$ we have $T\mathrm{Ext}_2^{s,t} = \{0\}$ for $t - s = 2n$ and thus we only have to deal with primes $p$ with $p(p-1) > n+1$. Since $n \geq 3$ this implies $p$ is odd. The proof is finished by using the well known result that $T \oplus_{t-s=2n} \mathrm{Ext}_A^{s,t}(\mathbb{Z}/p, \mathbb{Z}/p) = 0$ for odd primes with $(p-1)p > n+1$ (cf. [Na]). Starting from this, one proves inductively that if $M = \oplus M^k$ is a $(-1)$-connected graded $A$-module with $M^k$ finitely generated and $M^{2k+1} = 0$ for all $k \leq n$, then $T \oplus_{t-s=2n} \mathrm{Ext}_A^{s,t}(M; \mathbb{Z}/p) = 0$. Since $H^\star(S; \mathbb{Z}/p)$ fulfills this condition the proof of Proposition 12 is finished.

## 9. Analysis of certain obstructions in $l_0(\{e\})$

In this section we study a special class of 7-manifolds with very simple normal 2-type. We will show that under appropriate assumptions the $l_8$ $(= l_0)$-obstruction for a $B$-bordism between such manifolds can be controlled by characteristic numbers. The special class of manifolds is motivated by the classification of certain homogeneous spaces which fall into this class [Kr-St1], [Kr-St2], [Kr-St3]. Let $M$ be a 1-connected compact 7-dimensional manifold with $H_2(M; \mathbb{Z})$ torsion free of rank $r$. The normal 2-type $B$ of such manifolds is

$$\xi : \underbrace{\mathbb{C}P^\infty \times \ldots \times \mathbb{C}P^\infty}_{r\text{-copies}} \times B\operatorname{Spin} \xrightarrow{H(w_2) \oplus p} BO,$$

where $p : B\operatorname{Spin} \longrightarrow BO$ is the canonical projection. The map $H(w_2)$ is trivial, if $w_2(M) = 0$ or the classifying map for the Hopf bundle over one copy of $\mathbb{C}P^\infty$, if $w_2(M) \neq 0$. Now $H^4(B\operatorname{Spin}; \mathbb{Z}) \cong \mathbb{Z}$ is generated by a class denoted $\frac{p_1}{2}$. Thus for a spin vector bundle the characteristic class $\frac{p_1}{2}$ can be defined as the pullback of this class.

THEOREM 6. *Let $M_0$ and $M_1$ be 1-connected 7-dimensional, compact manifolds, either both* Spin *or both non-*Spin *such that $H_2(M_i; \mathbb{Z})$ is torsion free of rank $r$ and $H^4(M_0; \mathbb{Z}) \cong H^4(M_1; \mathbb{Z})$ is finite and generated by products of classes in $H^2(M_i; \mathbb{Z})$ and $\frac{p_1}{2}(M_i)$, if $M_i$ is* Spin *or $\frac{p_1}{2}(TM_i \oplus L)$, if $M_i$ is not* Spin *and $L$ is a complex line bundle with $w_2(L) = w_2(M_i)$.*

*Then $M_0$ is diffeomorphic to $M_1$ if and only if there exist normal $B$-smoothings of $M_i$ and a $B$-bordism $(W, \overline{\nu})$ with*

(i) sign $W = 0$.
(ii) $\langle \overline{\nu}^\star x \cup \overline{\nu}^\star y, [W, \partial W] \rangle = 0$ *for all $x, y \in H^4(B; \mathbb{Q})$.*

The second condition is to be understood as follows. As $H^3(\partial W; \mathbb{Q}) = 0 = H^4(\partial W; \mathbb{Q})$ there is an isomorphism $H^4(W, \partial W; \mathbb{Q}) \xrightarrow{\cong} H^4(W; \mathbb{Q})$. So regard $\overline{\nu}^\star x$ and $\overline{\nu}^\star y$ as elements in $H^4(W, \partial W; \mathbb{Q})$ before taking the cup product and evaluating on the fundamental class $[W, \partial W]$. As $H^4(B; \mathbb{Z})$ is generated by $\frac{p_1}{2}$ and products $z_i \cup z_j$, where the $z_i$'s generate the second cohomology of $B$, the condition on $H^4(M_i; \mathbb{Z})$ means that $H^4(M_i; \mathbb{Z})$ is finite and that $\overline{\nu}^\star : H^4(B; \mathbb{Z}) \to H^4(M_i; \mathbb{Z})$ is surjective. This result implies Theorem G. A generalization of it was proved in [Be].

*Proof.* By Proposition 4 we can assume that $W \xrightarrow{\overline{\nu}} B$ is a 4-equivalence. Then by Theorem 3 the surgery obstruction $\theta(W, \overline{\nu})$ for transforming $W$ into an $h$-cobordism is given by

$$\theta(W, \bar{\nu}) := (H_4(W, M_0; \mathbb{Z}) \xleftarrow{f} \operatorname{im}(d : \pi_5(B, W)$$
$$\to \pi_4(W)) \xrightarrow{g} H_4(W, M_1; \mathbb{Z}), \lambda, \mu).$$

Since $\mu$ is, in our situation, determined by $\lambda$ we omit it. We are going to show that $\theta(W, \bar{\nu})$ is elementary. Since the intersection form can be better treated in cohomology we will translate $\theta(W, \bar{\nu})$ to cohomology:

$$\theta(W, \bar{\nu}) := (H^4(W, M_1; \mathbb{Z}) \leftarrow KH^4(W, \partial W; \mathbb{Z}) \to H^4(W, M_0; \mathbb{Z}), \lambda),$$

where $KH^4(W, \partial W; \mathbb{Z}) = \operatorname{Ker} \rho : H^4(W, \partial W; \mathbb{Z}) \cong H_4(W; \mathbb{Z}) \to H_4(B; \mathbb{Z})$ and $\lambda$ is the cup-product pairing between $H^4(W, M_1)$ and $H^4(W, M_0)$. From the long exact sequences for the pair $(W, M_i)$ we see that $H^4(W, M_i; \mathbb{Z})$ can be considered as a kernel of $H^4(W; \mathbb{Z}) \longrightarrow H^4(M_i; \mathbb{Z})$.

Now we introduce some notation: $V := H^4(W, \partial W; \mathbb{Z})$, $A := H_4(B; \mathbb{Z})$, $\rho : V \to A$ the map above. We identify $H^4(M_0; \mathbb{Z})$ by some isomorphism with $H^4(M_1; \mathbb{Z})$ and denote this finite abelian group by $H$. We identify $H^4(W; \mathbb{Z})$ with $V^\star$ via Kronecker isomorphism and Poincaré duality and denote the adjoint of the intersection form by $S : V \to V^\star$. The cohomology sequence of the pair $(W, \partial W)$ translates into a short exact sequence:

$$0 \to V \xrightarrow{S} V^\star \xrightarrow{j} H \oplus H \to 0.$$

We denote the projection to the $i$-th factor $(i = 0$ or $1)$ by $p_i : H \oplus H \to H$. Then our obstruction $\theta(W, \bar{\nu})$ translates to

$$\theta(W, \bar{\nu}) = (\operatorname{Ker} p_1\, j \xleftarrow{f} \operatorname{Ker} \rho \xrightarrow{g} \operatorname{Ker} p_0 j, \lambda) \in l_0(\{e\}),$$

where the maps $f$ and $g$ are the restriction of $S$ to $\operatorname{Ker} \rho$ ($\operatorname{Ker} \rho$ maps under $S$ injectively to $\operatorname{Ker} p_0 j \cap \operatorname{Ker} p_1 j$) and $\lambda$ is induced by $S$. We note that the dual sequence of $0 \to \operatorname{Ker} \rho \xrightarrow{i} V \xrightarrow{\rho} A \to 0$ is again exact (since the groups are free):

$$0 \longleftarrow (\operatorname{Ker} \rho)^\star \xleftarrow{i^\star} V^\star \xleftarrow{\rho^\star} A^\star \longleftarrow 0.$$

The following algebraic data can be derived from our topological assumptions:

a) $S$ is symmetric and nondegenerate and sign $S = 0$.
b) If $S(v_i) = \rho^*(\alpha_i)$ for $v_i \in V, \alpha_i \in A^*$: $(i = 1, 2)$, then $S(v_1)(v_2) = 0$.
c) For all $\varphi \in V^*$ satisfying $p_1 j(\varphi) \neq 0$ and $p_0 j(\varphi) = 0$ we have $\rho(S^{-1}(r\varphi)) \neq 0$ for all positive multiples $r\varphi$ of $\varphi$ lying in the image of $S$ (and the same holds if we interchange the indices 0 and 1).

Now a) follows from sign $W = 0$, b) is a consequence of assumption (ii), and c) follows from the unimodularity of the linking forms on $M_i$:

Let $L$ denote the cohomology linking form of $M_1$. For a contradiction we assume $\varphi \in H^4(W; \mathbb{Z}) = V^*$ and a nontrivial $\psi \in H^4(M_1; \mathbb{Z}) = H \oplus \{0\}$ with

$$j\varphi = (\psi, 0) \in H \oplus H = H^4(M_1; \mathbb{Z}) \oplus H^4(M_0; \mathbb{Z}) = H^4(\partial W; \mathbb{Z})$$

such that for some $r \in \mathbb{N}$ we have $r\varphi \in \mathrm{Im}\,S$ and $\rho S^{-1}(r\varphi) = 0$. We will show then that $\psi = 0$, contradicting our assumption. Since the linking form $L$ on $H^4(M_1)$ is unimodular this is equivalent to showing that $L(\alpha, \psi) = 0$ for all $\alpha \in H^4(M_1)$. Since $\bar{\nu}_1^* : H^4(B; \mathbb{Z}) \longrightarrow H^4(M_1; \mathbb{Z})$ is surjective by assumption, we have to check $L(\bar{\nu}_0^*\eta, \psi) = 0$ for all $\eta \in H^4(B; \mathbb{Z})$. In the situation above the relation between the intersection form on $W$ and the linking form on $\partial W$ implies:

$$L(\bar{\nu}_0^*\eta, \psi) = \pm\frac{1}{r}\,\langle \eta, \rho S^{-1}(r\varphi) \rangle.$$

The latter expression vanishes since $\rho S^{-1}(r\varphi) = 0$.

The following proposition finishes the proof of Theorem 6.

PROPOSITION 13.  *Let $V$, $A$, $\rho : V \to A$ and $S$ be as introduced above. Then, if the assumptions* a)–c) *are fulfilled*

$$\theta(W, \bar{\nu}) = (\mathrm{Ker}\,p_1 j \xleftarrow{\ f\ } \mathrm{Ker}\,\rho \xrightarrow{\ g\ } \mathrm{Ker}\,p_0 j, \lambda) \in l_0(\{e\})$$

*is elementary.*

*Proof.* We denote the restriction of the symmetric bilinear form $S$ on $V$ to $\mathrm{Ker}\,\rho$ by $S_K$. The adjoint of $S_K$ is given by the composition $\mathrm{Ker}\,\rho \xrightarrow{i} V \xrightarrow{\rho} V^* \xrightarrow{i^*} (\mathrm{Ker}\,\rho)^*$.

Using property b) one shows:

i) For $v \in V, \alpha \in A^*$, $S(v) = \rho^*(\alpha) => \rho(v) = 0$.

From this and the definition of the radical one has:

ii) $S^{-1}(n \cdot \mathrm{Im}\,\rho^*) \subseteq \mathrm{rad}\,(S_K)$, where $n$ is the exponent of $H$ and rad is the radical and $\mathrm{rank}\,(\mathrm{rad}\,(S_K)) = \mathrm{rank}\,A^* = \mathrm{rank}\,A$.

Finally we will show:

iii) Cokernel $S_K$ is torsion-free; hence the form $\tilde{S}_K$ on $\mathrm{Ker}\,\rho/\mathrm{rad}\,(S_K)$ induced by $S_K$ is unimodular and its signature vanishes.

We will prove iii) at the end and finish the proof of Proposition 13 using i)–iii). By iii) there exists $U \subseteq \mathrm{Ker}\,\rho$ such that $U \cap \mathrm{rad}\,(S_K) = 0$ and $U = U^\perp$ is a direct summand of $\mathrm{Ker}\,\rho/\mathrm{rad}\,(S_K)$ of half rank. We show that $U \oplus \mathrm{rad}\,(S_K)$ maps under $f$ and $g$ in $\mathrm{Ker}\,(p_i j)$ to direct summands $B_i$ of half rank, proving the proposition.

It is clear that they have the right rank since $\mathrm{rank}\,\mathrm{Ker}\,(p_i j) = \mathrm{rank}\,(V^*) = \mathrm{rank}\,V = \mathrm{rank}\,(\mathrm{Ker}\,\rho) + \mathrm{rank}\,A = 2 \cdot \mathrm{rank}\,(U) + \mathrm{rank}\,(\mathrm{rad}\,(S_K)) + \mathrm{rank}\,A =$ $2(\text{rank}\,(U) + \text{rank}\,(\text{rad}\,S_K))$. To show that they are direct summands we first note that $U \oplus \text{rad}\,(S_K)$ is a direct summand in $\text{Ker}\,\rho$ and thus in $V$. This implies that if $x \in \text{Ker}\,p_0 j$ represents a nontrivial torsion element in $\text{Ker}\,p_0 j / f(U \oplus \text{rad}\,(S_K))$, then $x$ cannot be in the image of $S$ or equivalently $p_1 j(x) \neq 0$. But then by c) we obtain a contradiction since, as $x$ represents a torsion element, some multiple $rx = S(y)$ for some $y$ in $U \oplus \text{rad}\,S_k \subset \text{Ker}\,\rho$. The same argument holds if we consider $\text{Ker}\,p_1 j$.

We finish the argument by showing iii). Denote the inclusion from $\text{Ker}\,\rho$ to $V$ by $i$. Assume cokernel $S_K$ has torsion. Then there would exist $\varphi \in V^*$ with $i^*(\varphi) \notin \text{im}(S_K)$, but $i^*(r\varphi) = S(y)$ for some $r \in \mathbb{N}$ and $y \in \text{Ker}\,\rho$. Consider the two cases $j(\varphi) \in \text{im}(j\rho^*)$ and $j(\varphi) \notin \text{im}\,(j\rho^*)$. We will show that both lead to a contradiction.

If $j(\varphi) \in \text{im}(j\rho^*)$, there is an $\alpha \in A^*$ with $j(\varphi) = j\rho^*(\alpha)$. Define $\varphi' := \varphi - \rho^*(\alpha)$. Then for $\varphi'$ we have $i^*(\varphi) = i^*(\varphi')$ and $j(\varphi') = 0$. The latter implies that $\varphi'$ has a pre-image $v \in V$ under $S$. On the other hand $S(r \cdot v - i(y)) \in \text{im}\,\rho^*$. By i) we get $0 = \rho(r \cdot v - i(y)) = r\rho(v)$. As $A$ is torsion free, it follows $v \in \text{Ker}\,\rho$, a contradiction.

If $j\varphi \notin \text{im}(j\rho^*)$ choose $\alpha \in A^*$ such that $p_0 j \rho^*(\alpha) = p_0 j(\varphi)$ and define $\varphi' := \varphi - \rho^*(\alpha)$. Then we have $i^*\varphi' = i^*\varphi$ and $p_0 j(\varphi') = p_0 j(\varphi) - p_0 j\,\rho^*(\alpha) = 0$. The assumption $j(\varphi) \notin \text{Im}(j\rho^*)$ implies $p_1 j(\varphi') \neq 0$. Next we will show that $\rho S^{-1}(n \cdot r\varphi') = 0$, giving a contradiction to c). As above, we conclude from i) that $S^{-1}(nr\varphi' - i(n \cdot y)) \in \text{Ker}\,\rho$ and hence $\rho S^{-1}(n \cdot r\varphi') = 0$.

To show that $\text{sign}(S_K) = 0$, choose $X \subseteq \text{Ker}\,\rho$ such that $\text{Ker}\,\rho = \text{rad}\,(S_K) \oplus X$. This is possible because $\text{Ker}\,\rho/\text{rad}\,(S_K)$ is free. As $S_K|_{X \times X}$ is unimodular, we have $V = X + X^\perp$. We can choose $Y \subset X^\perp$ such that $\rho|_Y : Y \longrightarrow A$ is an isomorphism, since $A$ is free. Starting from the decomposition $V = X \oplus \text{rad}\,(S_K) \oplus Y$ we note that $\text{sign}(S|_{\text{rad}\,(S_K) \oplus Y})$ is zero, because $\text{rad}\,(S_K)$ has the same rank as $Y$ and the form vanishes on $\text{rad}\,(S_K)$. On the other hand $X$ is orthogonal to $\text{rad}\,(S_K) \oplus Y$ and thus

$$0 = \text{sign}(S) = \text{sign}(S|_X) + \text{sign}(S|_{\text{rad}\,(S_K) \oplus Y}) = \text{sign}(S|_X) = \text{sign}(S_K). \qquad \square$$

Fachbereich Mathematik, Universität Mainz, 55099 Mainz
and
Mathematisches Forschungsinstitut Oberwalfach, 77709 Oberwolfach,
Federal Republic of Germany
*E-mail address:* kreck@mfo.de